\documentclass[journal]{IEEEtran}
\ifCLASSINFOpdf
\else
\fi

\usepackage{amsbsy}
\usepackage{amsfonts}
\usepackage{amsmath}
\usepackage{graphicx}
\usepackage{multirow}
\usepackage{algorithm}
\usepackage[noend]{algpseudocode}
\usepackage{cite}
\usepackage{epstopdf}
\usepackage{color}
\usepackage{float}
\usepackage{textcomp}

\makeatletter
\let\myorg@bibitem\bibitem
\def\bibitem#1#2\par{%
  \@ifundefined{bibitem@#1}{%
    \myorg@bibitem{#1}#2\par
  }{%
    \begingroup
      \color{\csname bibitem@#1\endcsname}%
      \myorg@bibitem{#1}#2\par
    \endgroup
  }%
}
\makeatother

\newcommand{\tabincell}[2]{\begin{tabular}{@{}#1@{}}#2\end{tabular}}

\begin{document}
%
\title{An  SQP Method Combined with Gradient Sampling for Small-Signal Stability Constrained OPF }
%
%
%

\author{Peijie~Li,~\IEEEmembership{Member,~IEEE,}
        Junjian~Qi,~\IEEEmembership{Member,~IEEE,}
        Jianhui~Wang,~\IEEEmembership{Senior~Member,~IEEE,}
        Hua~Wei,
        Xiaoqing~Bai,
        and Feng~Qiu~\IEEEmembership{}

\thanks{This work was supported in part by National Natural Science
Foundation of China under Grant 51407036, 51367004, the National Basic Research Program of China (973 Program) under Grant 2013CB228205 and the U.S. Department of Energy (DOE)'s Office of Electricity Delivery and Energy Reliability.

P. Li is with the College of Electrical  Engineering, Guangxi University, Nanning, 530004 China. He is visiting Argonne National Laboratory, Lemont, IL 60439 USA (e-mail: lipeijie@gxu.edu.cn).}
\thanks{J. Qi, J. Wang, and F. Qiu are with Energy Systems Division, Argonne National Laboratory, Lemont, IL 60439 USA (e-mails: jqi@anl.gov; jianhui.wang@anl.gov; fqiu@anl.gov).}
\thanks{H. Wei and X. Bai are with the College of Electrical  Engineering, Guangxi University, Nanning, 530004 China.}
\thanks{}
\vspace*{-0.5cm}
}

\markboth{Preprint of DOI 10.1109/TPWRS.2016.2598266, IEEE Transactions on Power Systems.}{stuff}

\maketitle

\begin{abstract}
Small-Signal Stability Constrained Optimal Power Flow (SSSC-OPF) can provide additional stability measures and control strategies to guarantee  the system to be small-signal stable. However, due to the nonsmooth property of the spectral abscissa function, existing  algorithms solving SSSC-OPF cannot guarantee convergence. To tackle this computational challenge  of SSSC-OPF, we propose a Sequential Quadratic Programming (SQP) method combined with Gradient Sampling (GS) for SSSC-OPF. At each iteration of the proposed SQP, the gradient of the spectral abscissa function is randomly sampled at the current iterate and additional nearby points to make the search direction computation effective in nonsmooth regions. The method can guarantee SSSC-OPF is globally and efficiently convergent to stationary points with probability one. The effectiveness of the proposed method is tested and validated on WSCC 3-machine 9-bus system, New England 10-machine 39-bus system, and IEEE 54-machine 118-bus system.
\end{abstract}

\begin{IEEEkeywords}
small-signal stability, AC optimal power flow (ACOPF), gradient sampling, nonsmooth optimization,  sequential quadratic programming,  spectral abscissa.
\end{IEEEkeywords}

%
\IEEEpeerreviewmaketitle

\section{Introduction}
%
%
%
%

\IEEEPARstart{I}{n} the large interconnected power systems, the small-signal stability problem with the oscillation behavior can greatly threaten the system security \cite{power_system_oscillation}. 
The increased penetration of renewable energy sources can further deteriorate the small-signal stability \cite{6737307,7053959,7337588}, which 
is mainly driven by:
\begin{enumerate}
\item Converter control-based generators including variable-speed doubly fed asynchronous generators and photovoltaic (PV) plants are replacing the conventional generators.  
These solid-state inverters, however, usually do not contribute to the system inertia \cite{6737307}.
\item The fluctuating nature of the renewable resources such as wind and solar may cause rapid changes in future generation patterns, leading to rapid fluctuations of the power system's operating point \cite{7053959}.
\item In some countries, such as China, due to the differences between the geographical distribution of load center and renewable energy, 
new Ultra-high-voltage Alternating Current (UHVAC) transmission lines are built to transfer the clean energy. 
The long distance power transmission over these UHVAC lines can result in oscillations \cite{7337588}.  
\end{enumerate}

Although the damping controllers can  enhance small-signal stability, they cannot always guarantee the system to be small-signal stable \cite{1266608}. 
By  contrast,  re-dispatch can provide additional stability measures and control strategies to make the system to be small-signal stable. 
Most commercial software, such as PSS\textsuperscript{\textregistered}E and NEPLAN, only uses participation factors \cite{sauer1998power,kundurpower} to perform offline re-dispatch study. 
Nonetheless, participation factors  can neither determine whether the generation for each generator should be increased or decreased, nor tell how much generation should be dispatched \cite{chinazhou}.
On the contrary, the Small-Signal  Stability  Constrained Optimal Power Flow (SSSC-OPF) is a commendable model that can give the complete re-dispatch information to guarantee the small-signal stability  while considering the economic objective and  technical  constraints. However, solving the SSSC-OPF problem can be very challenging because of the nonsmooth nature of the small-signal stability constraint. 
Existing methods for solving SSSC-OPF include:
\begin{itemize}
\item Numerical eigenvalue sensitivity based Interior Point Method (IPM): Numerical eigenvalue sensitivities have been used to improve the power transfer capability constrained by small-signal stability \cite{1266608}. In \cite{5593195}, an SSSC-OPF is formulated to achieve an appropriate security level under stressed loading conditions, 
in which the small-signal stability constraint is replaced by first-order Taylor series expansion and the gradients of the real part of the critical eigenvalue are computed by numerical eigenvalue sensitivities. 
This method tries to make the small-signal stability constraint feasible by successively solving SSSC-OPF with IPM, but it may compromise on the economic cost due to losing the high-order terms of the small-signal constraint.

\item  Approximate singular value sensitivity based IPM: An SSSC-OPF is used to tune the oscillation controls in electricity markets \cite{4349096},  
in which the minimum singular value of the modified full Jacobian matrix is proposed as an stability index. 
However, the gradient of this index is derived approximately by first-order Taylor series expansion and the Hessian is numerically evaluated by perturbing the gradients. 

\item Closed-form eigenvalue sensitivity based IPM: An expected-security-cost OPF  with the small-signal stability constraint is presented in \cite{1717577}, 
in which a more computational efficient closed-form formula is used to calculate the eigenvalue sensitivities \cite{alvarado1999avoiding}, \cite{852145}. 
Nevertheless, the second-order eigenvalue sensitivities essential for forming the Hessian of the small-signal stability constraints for IPM have to be derived and calculated, which can be very time-consuming. 

\item Nonlinear semi-definite programming (NLSDP) method: A nonlinear semi-definite programming model is proposed to formulate the  spectral abscissa constraint as a semi-definite constraint indirectly which is further transformed to some smooth nonlinear constraints \cite{6531914}. An explicit and equivalent small-signal stability constraint is obtained based on Lyapunov theorem. However, the dense subsidiary semi-definite matrix variables may make the model computationally prohibitive for large systems.
\end{itemize}

IPM methods cannot guarantee convergence of SSC-OPF, even locally, because of the nonsmooth nature of the small-signal stability constraint function. 
During  iterations they may suffer from oscillations between critical modes for different generator rescheduling patterns.  
Recently, a  Sequential Quadratic Programming (SQP) method combined with Gradient Sampling (GS), called SQP-GS, 
is proposed for the nonsmooth constrained optimization \cite{curtis2012sequential}. 
Based on the SQP-GS, this paper  proposes an optimization method to solve the SSSC-OPF problem with global convergence. 
 The closed-form  eigenvalue sensitivity is used to calculate the gradient of the spectral abscissa. 
 Since the SQP only needs the gradients of all functions in the model, the  second-order eigenvalue sensitivities are not needed. 
 Moreover, owing to  only one nonsmooth function in SSSC-OPF, GS can be performed only for this nonsmooth function, which would be beneficial to  improve computation efficiency significantly. 
 
The remainder of this paper is organized as follows. Section \ref{sssm} describes the small-signal stability model. In Section \ref{sss}, the property of the small-signal stability constraint and the calculation of the eigenvalue sensitivities are discussed. Section \ref{sqp-gs} introduces the SQP-GS method for nonsmooth constrained optimization. Then we discuss how to solve the SSSC-OPF problem by the SQP-GS method in Section \ref{sss-opf} . The proposed method is tested and validated on three systems in Section \ref{result}. Finally the conclusion is drawn in Section \ref{conclusion}.
\section{Small-Signal Stability Model } \label{sssm}
\subsection{General Model of Small-Signal Stability}
   The behavior of a dynamical system can be described by differential and algebraic equations in the following form
    \begin{align}
    \label {da:1}   \dot{\boldsymbol s} &= \boldsymbol f_d(\boldsymbol s,\boldsymbol u,\boldsymbol y)\\
    \label {da:2}    \boldsymbol 0 &= \boldsymbol f_a(\boldsymbol s,\boldsymbol y),
     \end{align}
   where $\boldsymbol s$ is the vector of state variables, $\boldsymbol u$ is the vector of inputs, and $\boldsymbol y$ is the non-state variables.
   
   In small-signal stability analysis, the linearized form of (\ref{da:1})--(\ref{da:2}) is often used: \\
   \begin{equation}
    \label {lin:1} 
    \begin{bmatrix}
       \Delta \dot{\boldsymbol s} \\
      \boldsymbol 0
      \end{bmatrix}
      =
        \begin{bmatrix}
       \tilde{\boldsymbol A}   &\tilde{\boldsymbol B} \\
       \tilde{\boldsymbol C}   &\tilde{\boldsymbol D}
      \end{bmatrix} 
    \begin{bmatrix}
       \Delta \boldsymbol s   \\
       \Delta \boldsymbol y  
      \end{bmatrix}
      + 
          \begin{bmatrix}
       \boldsymbol E_1 \\
       \boldsymbol 0 
      \end{bmatrix}
       \Delta \boldsymbol u.
\end{equation}   
  Eliminating $\Delta \boldsymbol y$ we can get 
      \begin{eqnarray}  \label {eq:diff}
    \boldmath \Delta \dot{\boldsymbol s}=\boldsymbol A\Delta \boldsymbol s +  \boldsymbol E_1 \Delta \boldsymbol u,
   \end{eqnarray}
   where $\boldsymbol A= \tilde{\boldsymbol A}-\tilde{\boldsymbol B}\tilde{\boldsymbol D}^{-1}\tilde{\boldsymbol C}$  is commonly referred to as the state matrix.

\subsection{Differential and Algebraic  Equations of Power Systems}

    \subsubsection{Generator Model}

    The two-axis synchronous generator model \cite{sauer1998power} that has been widely used in small-signal stability analysis is considered in this paper. This model is more realistic than the classical model used in \cite{oldsen}, which is appropriate only for the most basic studies.
    The differential equations for generator $i \in S_\textrm{G}$ can be written as
    \begin{align}
    \label{gen:delta} &\frac{d\delta_i}{dt} = \omega_i - \omega_\textrm{s} \\
    \label{gen:omega} &\frac{d\omega_i}{dt} =  \frac{1}{M_i} \Big(T_{\textrm{M}i}- (E_{\textrm{q}i}^{'}- X_{\textrm{d}i}^{'} I_{\textrm{d}i}) I_{\textrm{q}i}\notag\\
    &\qquad\qquad -( E_{\textrm{d}i}^{'}+ X_{\textrm{q}i}^{'} I_{\textrm{q}i}) I_{\textrm{d}i}-D_i( \omega_i- \omega_s)\Big)  \\
        \label{gen:Eq} &\frac{dE_{\textrm{q}i}^{'}}{dt} = \frac{1}{T_{\textrm{d0}i}^{'}} \Big( -E_{\textrm{q}i}^{'} 
    - (X_{\textrm{d}i}- X_{\textrm{d}i}^{'})I_{\textrm{d}i}+E_{\textrm{fd}i} \Big)\\
    \label{gen:Ed} &\frac{dE_{\textrm{d}i}^{'}}{dt} = \frac{1}{T_{\textrm{q0}i}^{'}} \Big( -E_{\textrm{d}i}^{'} 
    +(X_{\textrm{q}i}-X_{\textrm{q}i}^{'})I_{\textrm{q}i}\Big),
    \end{align}
     where
    $S_\textrm{G}$ is the set of generators, $\delta_i$ is rotor angle, $\omega_i$  is rotor speed, $\omega_\textrm{s}$ is rated rotor speed, 
     $E_{\textrm{d}i}^{'}$ and $E_{\textrm{q}i}^{'}$ are, respectively, the d-axis and q-axis components of the internal voltage, 
     $I_{\textrm{d}i}$ and $I_{\textrm{q}i}$ are d-axis and q-axis components of the internal current,
     $T_{\textrm{M}i}$  is mechanical power output, $M_i$ is inertia constant, $D_i$ is damping torque coefficient, 
     $E_{\textrm{fd}i}$  is excitation output voltage,
     $X_{\textrm{d}i}$ and  $X_{\textrm{q}i}$  are synchronous reactance,
     $X_{\textrm{d}i}^{'}$ and  $X_{\textrm{q}i}^{'}$  are transient reactance, and 
     $T_{\textrm{d0}i}^{'}$ and $T_{\textrm{q0}i}^{'}$ are open-circuit time constant, respectively, at d and q axes.

   Besides, the stator algebraic equations for generator $i \in S_\textrm{G}$  can be written in polar form as \cite{sauer1998power}
      \begin{align}
        \label{gen_al:Ed} & E_{\textrm{d}i}^{'}-  V_i\sin(\delta_i- \theta_i)- R_{\textrm{s}i} I_{\textrm{d}i}+ X_{\textrm{q}i}^{'} I_{\textrm{q}i} = 0 \\
      \label{gen_al:Eq} & E_{\textrm{q}i}^{'}- V_i\cos(\delta_i- \theta_i)- R_{\textrm{s}i} I_{\textrm{q}i}- X_{\textrm{d}i}^{'} I_{\textrm{d}i} = 0,
      \end{align}
     where $V_{i}$  is bus voltage magnitude, $\theta_i$ is bus voltage phase angle, and $R_{\textrm{s}i}$ is the armature resistance.
    
     \subsubsection{Exciter Model}
     In this paper we use the IEEE Type DC-1 exciter for each generator, which can be expressed in differential equations for generator  $i \in S_\textrm{G}$ as \cite{sauer1998power}
     \begin{align}
      \label{Exc:efd} &\frac{dE_{\textrm{fd}i}}{dt} = \frac{1}{T_{\textrm{E}i}}\Big(  -\big(K_{\textrm{E}i}+S_E(E_{\textrm{fd}i})\big) E_{\textrm{fd}i}+V_{\textrm{R}i} \Big)\\
      \label{Exc:vr} &\frac{dV_{\textrm{R}i}}{dt} = \frac{1}{T_{\textrm{A}i}}\Big( -V_{\textrm{R}i}+K_{\textrm{A}i} R_{\textrm{F}i}-\frac{K_{\textrm{A}i}K_{\textrm{F}i}}{T_{\textrm{F}i}}E_{\textrm{fd}i} \notag\\
      &\qquad \qquad\qquad\qquad\qquad\quad\; +K_{\textrm{A}i}(V_{\text{ref}i}-V_i)\Big)\\
      \label{Exc:rf} &\frac{dR_{\textrm{F}i}}{dt} = 
      \frac{1}{T_{\textrm{F}i}} \Big( -R_{\textrm{F}i}+\frac{K_{\textrm{F}i}}{T_{\textrm{F}i}}E_{\textrm{fd}i}\Big),
     \end{align}
    where
        $V_{\textrm{R}i}$ is voltage regulator output,
        $R_{\textrm{F}i}$ is exciter rate feedback,
        $S_E( E_{\textrm{fd}i}) =   A_{\textrm{e}i}e^{B_{\textrm{e}i}E_{\textrm{fd}i}}$ is the field saturation function with coefficients $A_{\textrm{e}i}$ and $B_{\textrm{e}i}$, 
        $K_{\textrm{E}i}$ is exciter gain,
        $K_{\textrm{A}i}$ is voltage regulator gain,
        $K_{\textrm{F}i}$ is rate feedback gain,
        $T_{\textrm{E}i}$ is exciter time constant,
        $T_{\textrm{A}i}$ is voltage regulator time constant,
        $T_{\textrm{F}i}$ is rate feedback time constant,
        and $V_{\text{ref}i}$ is the reference voltage.

     \subsubsection{Network  Model}
     
     The  network equations relate the real and reactive power injections at each bus to the voltage magnitudes and phase angles at the system buses. 
     In this paper the loads are modeled as constant power. Then for a generator bus $i \in S_G$  we have \cite{sauer1998power}
     \begin{align}
     \label{NW:pg} I_{\textrm{d}i}V_i\sin(\delta_i-\theta_i)+I_{\textrm{q}i}V_i\cos(\delta_i-\theta_i)+ P_{\textrm{L}i}  \notag\\
     -\sum_{j\in S_\textrm{B}}V_iV_jY_{ij}\cos(\theta_i-\theta_j-\alpha_{ij})=0 \\
     \label{NW:qg}  I_{\textrm{d}i}V_i\cos(\delta_i-\theta_i)-I_{\textrm{q}i}V_i\sin(\delta_i-\theta_i)+ Q_{\textrm{L}i}  \notag\\
     -\sum_{j\in S_\textrm{B}}V_iV_jY_{ij}\sin(\theta_i-\theta_j-\alpha_{ij})=0,
     \end{align}
     where
     $P_{\textrm{L}i}$ and $Q_{\textrm{L}i}$ are, respectively, the active and reactive load, and $Y_{ij}e^{j\alpha_{ij}}$ is the entry of the admittance matrix.

     For a non-generator bus $i \in S_\textrm{L}$, there are
     \begin{align}
     \label{NW:pd}&P_{\textrm{L}i}-\sum_{j\in S_\textrm{B}}V_iV_jY_{ij}\cos(\theta_i-\theta_j-\alpha_{ij})=0 \\
     \label{NW:qd}&Q_{\textrm{L}i}-\sum_{j\in S_\textrm{B}}V_iV_jY_{ij}\sin(\theta_i-\theta_j-\alpha_{ij})=0,
     \end{align}
     where $S_\textrm{L}$ is the set of non-generator buses and $S_\textrm{B} = S_\textrm{G}\cup S_\textrm{L}$ is the set of all of the buses.
     
    \subsection {Linearization of Dynamic System Model}
   
    Linearization of (\ref{gen:delta})--(\ref{NW:qd}) yields
         \begin{align}
          \label{MT:1} \Delta \dot{\boldsymbol s} = &\,\boldsymbol A_1 \Delta \boldsymbol s + \boldsymbol B_1 \Delta \boldsymbol I_\textrm{g} + \boldsymbol B_2 \Delta \boldsymbol V_\textrm{g} +\boldsymbol E_1 \Delta \boldsymbol u \\         
           \label{MT:2} 0 = &\,\boldsymbol C_1 \Delta \boldsymbol s + \boldsymbol D_1 \Delta \boldsymbol I_\textrm{g} + \boldsymbol D_2 \Delta \boldsymbol V_\textrm{g}  \\   
            \label{MT:3}0 = &\,\boldsymbol C_2 \Delta \boldsymbol s + \boldsymbol D_3 \Delta \boldsymbol I_\textrm{g} + \boldsymbol D_4 \Delta \boldsymbol V_\textrm{g} + \boldsymbol D_5 \Delta  \boldsymbol V_\textrm{l}\\
            \label{MT:4}0 = &\,\boldsymbol D_6 \Delta \boldsymbol V_\textrm{g} + \boldsymbol D_7 \Delta  \boldsymbol V_\textrm{l},
         \end{align}
    where (\ref{MT:1}) is obtained by linearizing the differential equations (\ref{gen:delta})--(\ref{gen:Ed}) and  (\ref{Exc:efd})--(\ref{Exc:rf}), 
    (\ref{MT:2}) comes from the stator algebraic equations (\ref{gen_al:Ed})--(\ref{gen_al:Eq}), 
    (\ref{MT:3}) comes from the network equations (\ref{NW:pg})--(\ref{NW:qg}), 
     and (\ref{MT:4}) is from the network equations (\ref{NW:pd})--(\ref{NW:qd});
     $\boldsymbol A_1$, $\boldsymbol B_1$, $\boldsymbol B_2$, $\boldsymbol E_1$, $\boldsymbol C_1$, $\boldsymbol D_1$, $\boldsymbol D_2$, $\boldsymbol C_2$, $\boldsymbol D_3$ are block diagonal matrices, $\boldsymbol D_4$, $\boldsymbol D_5$, $\boldsymbol D_6$, $\boldsymbol D_7$ are full matrices, and    
     \begin{align}
           \boldsymbol s_i = \, &\,[\delta_i\;\; \omega_i\;\; E_{\textrm{q}i}^{'}\;\; E_{\textrm{d}i}^{'}\;\; E_{\textrm{fd}i}\;\; V_{\textrm{R}i}\;\; R_{\textrm{F}i}]^\top\quad &&i\in S_G\notag\\
     \boldsymbol I_{\textrm{g}i} =  &\,[I_{\textrm{d}i}\;\; I_{\textrm{q}i}]^\top  &&i\in S_G\notag\\  
     \boldsymbol V_{\textrm{g}i} = \, &\,[\theta_i\;\; V_i ]^\top \quad &&i \in S_G\notag\\
      \boldsymbol V_{\textrm{l}i} = \, &\,[\theta_i \;\; V_i ]^\top \quad &&i \in S_L\notag\\
           \boldsymbol u_i = \, &\,[ T_{\textrm{M}i} \;\; V_{\text{ref}i}]^\top \quad &&i \in S_G \notag.
     \end{align}
     We can rewrite (\ref{MT:1})--(\ref{MT:4}) in the compact form (\ref{lin:1}), and
      \begin{align}
            \boldsymbol y = \begin{bmatrix} \boldsymbol I_\textrm{g}^\top \;\; \boldsymbol V_\textrm{g}^\top \;\;  \boldsymbol V_\textrm{l}^\top \end{bmatrix}^\top \notag
      \end{align}
     \begin{align}
       \tilde{\boldsymbol A} = \boldsymbol A_1,\quad
       \tilde{\boldsymbol B} =  &\begin{bmatrix} \boldsymbol B_1 & \boldsymbol B_2   &\boldsymbol 0 \end{bmatrix}, \notag\\
       \tilde{\boldsymbol C} =  \begin{bmatrix} \boldsymbol C_1\\ \boldsymbol C_2\\  \boldsymbol  0 \end{bmatrix} ,\quad
        \tilde{\boldsymbol D} = &\begin{bmatrix}
         \boldsymbol D_1  &\boldsymbol D_2    &\boldsymbol 0  \\
         \boldsymbol D_3  &\boldsymbol D_4    &\boldsymbol D_5  \\
          \boldsymbol 0   &\boldsymbol D_6  &\boldsymbol D_7
         \end{bmatrix}.\notag
     \end{align}
 
\section{Small-Signal Stability Constraint and Eigenvalue Sensitivities} \label{sss}
   If the eigenvalues of the state matrix $\boldsymbol{A}$ 
   all have negative real parts, the power system is stable in small-signal stability sense.   
   An index $\eta$ called \textit{spectral abscissa}, which is the largest real part of the eigenvalues, is often used to describe the security margin: 
   \begin{equation}
       \eta(\boldsymbol A) =\max\{\text{Re}(\lambda):\lambda\in \boldsymbol \lambda(\boldsymbol A)\} = \text{Re}(\lambda_\eta),
   \end{equation}
 where $\boldsymbol\lambda(\boldsymbol A)$ represents all the eigenvalues of $\boldsymbol A$,  $\text{Re}(\lambda)$ is the real parts of the eigenvalues $\lambda$, and $\lambda_{\eta}$ is the eigenvalue with the largest real part. 
 The spectral abscissa determines the decay rate of the amplitude of the oscillation. 
 The smaller the spectral abscissa, the more stable the system is.

The state matrix $\boldsymbol A$ usually has complex eigenvalues due to its unsymmetrical characteristic. Generally, the spectral abscissa function is non-smooth\cite{lewis1996eigenvalue}. Fortunately, the spectral abscissa function has been proved to be locally Lipschitz and continuously differentiable on open dense subsets of $\mathbb{R}^n$ \cite{burke2001variational}, which means that it is continuously differentiable almost everywhere and its gradient can be easily obtained where it is defined by calculating first-order spectral abscissa sensitivities.

As for computing the spectral abscissa sensitivities, the numerical differentiation method is widely used, 
which performs eigenvalue analysis to get the spectral abscissa $\eta(\boldsymbol A)$ of the state matrix at the equilibrium point and then vary one variable $x_i$ by a small quantity $\varepsilon$ 
to get the perturbed state matrix $A_\varepsilon$ and its spectral abscissa $\eta( \boldsymbol A_\varepsilon)$. 
The spectral abscissa sensitivity with respect to $x_i$ can be approximated by 
\begin{equation}
\frac{\partial \eta}{\partial  x_i} \approx \frac{\eta( \boldsymbol A_\varepsilon)-\eta(\boldsymbol A)}{ \varepsilon}.
\end{equation}

The numerical differentiation method is easy to implement, but for large systems its calculation burden can be heavy due to the repetitive procedure. 
Also, the sensitivity with respect to the power of the slack bus cannot be obtained.
Alternatively, the spectral abscissa sensitivity can be obtained by closed-form formulas. 
Specifically, the $j$th eigenvalue sensitivity with respect to the $i$th variable $x_i$ can be written 
as \cite{alvarado1999avoiding,852145}:
\begin{equation}\label {eq:eigen_sen}
\frac{\partial\lambda_j}{\partial x_i}=\frac{\boldsymbol\psi_j\frac{\partial \boldsymbol A }{\partial  x_i}\boldsymbol\phi_j}{\boldsymbol\psi_j \boldsymbol \phi_j},
\end{equation}
where $\boldsymbol \psi_j$ and $\boldsymbol\phi_j$ are, respectively, the left and right eigenvectors of the eigenvalue $\lambda_j$, and 
\begin{equation}
\begin{aligned}
 \label {eigen:A}    \frac{\partial \boldsymbol A}{\partial  x_i} =  &\frac{\partial \tilde{\boldsymbol A} }{\partial x_i}- \frac{\partial \tilde{\boldsymbol B }}{\partial  x_i}\tilde{\boldsymbol D }^{-1}\tilde{\boldsymbol C }
 \\
 &+\tilde{ \boldsymbol B }\tilde {\boldsymbol D }^{-1}\frac{\partial \tilde{\boldsymbol D }}{\partial   x_i}\tilde {\boldsymbol D }^{-1}\tilde{ \boldsymbol C }-\tilde{ \boldsymbol B }\tilde{\boldsymbol D }^{-1}\frac{\partial \tilde{\boldsymbol C }}{\partial  x_i}.
 \end{aligned}
\end{equation}
Then the sensitivity of the spectral abscissa with respect to $x_i$ can be given by
\begin{equation}
 \label {eigen:sa}    \frac{\partial \eta}{\partial x_i} = \text{Re}\Big(\frac{\partial\lambda_{\eta}}{\partial  x_i}\Big).
\end{equation}

From (\ref{eigen:A}) it is seen that the derivation of the eigenvalue sensitivities for all of the variables requires considerable work since the elements of 
$\tilde{\boldsymbol A }$, $\tilde{\boldsymbol B }$, $\tilde{\boldsymbol C}$, and $\tilde{\boldsymbol D }$ can be different functions of several variables. 
However, the derivation is required 
only once. 
In this paper, we use the closed-form formulation to compute the gradient of the  spectral abscissa.

\section{A Sequential Quadratic Programming Algorithm Combined with Gradient Sampling} \label{sqp-gs}  

Sequential Quadratic Programming (SQP) has a long and rich history in solving smooth constrained optimization problems \cite{boggs1995sequential}. 
In each iterate of the traditional SQP algorithms, a quadratic programming (QP) subproblem is solved to obtain a search direction. 
However, the traditional SQP algorithms will fail for nonsmooth problems in theory and in practice.  
In 2005 an algorithm known as Gradient Sampling (GS) was developed for nonsmooth unconstrained optimization problems \cite{burke2005robust}. 
More recently, a SQP-GS method that combines the techniques of SQP and GS is developed \cite{curtis2012sequential}, 
which is proved to be able to globally convergent to stationary points with probability one when the objective and constraint functions are locally Lipschitz and continuously differentiable on open dense subsets of $\mathbb{R}^n$. SQP-GS has been shown to be a reliable method for many challenging nonsmooth problems, even when the objective function is not locally Lipschitz, 
in which case the convergence cannot be guaranteed though \cite{curtis2012sequential}. 

The GS algorithm is conceptually simple. Basically, it is a stabilized steepest descent algorithm \cite{nocedal2006numerical}. 
In each iteration, a descent direction is obtained by evaluating the gradient of the objective function at the current iterate and an additional set of nearby points and computing the vector in the convex hull of the gradients with the smallest norm. A standard line search is then used to obtain a lower point. 
The stabilization is controlled by the gradient sampling radius. 
As a natural extension to constrained optimization, the GS procedure samples the gradients of the constraint functions along with the objective function, 
thereby ensuring that good search directions are produced in nonsmooth regions.

Generally, the SQP-GS algorithm is developed to solve optimization problems in the following form:
\begin{equation}
\label {model:nonlin} 
    \begin{aligned}
    & \underset{\boldsymbol x}{\text{min}}
    & & f(\boldsymbol x) \\
    & \text{s.t.}
    & & \boldsymbol h(\boldsymbol x) = 0 \\
    & & & \underline{\boldsymbol g} \leq \boldsymbol g(\boldsymbol x) \leq \overline{\boldsymbol g},
    \end{aligned}
\end{equation}
where the objective function $f: \mathbb{R}^n \to \mathbb{R}$, the equality constraint functions $\boldsymbol h:\mathbb{R}^n \to \mathbb{R}^s$, and the inequality constraint functions $\boldsymbol g:\mathbb{R}^n \to \mathbb{R}^m$ are locally Lipschitz and continuously differentiable on open dense subsets of $\mathbb{R}^n$. 

    At the heart of SQP-GS is the following   QP used to compute a search direction in the $k$th iteration:
      \begin{align}  
 \label {model:1}  \hspace*{-0.2cm}&\underset{\boldsymbol d, z,\boldsymbol e,\overline{\boldsymbol r},\underline{\boldsymbol r}}{\text{min}} \;\;\;\rho z+ \sum_{j=1}^{m_g} (\overline{r}^j+\underline{r}^j)+\sum_{i=1}^{m_h} e^i+\frac{1}{2}\boldsymbol d_k^\top \boldsymbol H_k\boldsymbol d_k\\
    \label {model:2}  &\text{s.t.}\quad\; f(\boldsymbol x_k)+\nabla f(\boldsymbol x)^\top\boldsymbol d_k    \leq z \qquad \;\;\;\;\;\;\,   \forall\boldsymbol x \in  \mathcal{B}^{f}_{\epsilon,k} \\ 
         \label {model:3} &\quad\quad\;  h^i(\boldsymbol x_k)+\nabla  h^i(\boldsymbol x)^\top\boldsymbol d_k  \leq  e^i\quad  \quad \; \;\;\,\,  \forall \boldsymbol x \in  \mathcal{B}^{h^i}_{\epsilon,k}\\
         \label {model:4} &\quad\; -h^i(\boldsymbol  x_k)-\nabla   h^i(\boldsymbol  x)^\top\boldsymbol d_k \leq -e^i \;\;\;\;\quad  \forall \boldsymbol x \in  \mathcal{B}^{h^i}_{\epsilon,k}\\
        \label {model:5}  &\quad\quad\;   g^j(\boldsymbol  x_k)+\nabla   g^j(\boldsymbol x)^\top\boldsymbol d_k  \leq \overline{r}^j+\overline{g}^j\;\;\; \forall \boldsymbol x \in  \mathcal{B}^{g^j}_{\epsilon,k}\\
    	 \label {model:6}&\quad\; -g^j(\boldsymbol  x_k)-\nabla   g^j(\boldsymbol x)^\top\boldsymbol d_k \leq \underline{r}^j-\underline{g}^j \;\; \forall \boldsymbol x \in \mathcal{B}^{g^j}_{\epsilon,k}\\
    	\label {model:7} &\quad\quad\;  (\boldsymbol e,\overline{\boldsymbol  r},\underline{\boldsymbol  r})  \geq 0,
    \end{align}
where
$\rho$ is a penalty parameter, $\boldsymbol d_k$ is the search direction, $\boldsymbol H_k$ is the approximated  Hessian of the Lagrangian of (\ref{model:nonlin}), 
$m_h$ and $m_g$ are, respectively, the number of the equality and inequality constraints, $z$, $\boldsymbol e$, $\overline{\boldsymbol r}$, and $\underline{\boldsymbol r}$ are slack variables, and 
\begin{align}
 \label {sample:1}\mathcal{B}_{\epsilon,k}^{f} &:= \{\boldsymbol x_{k,0}^f,\boldsymbol x_{k,1}^f,\cdots,\boldsymbol x_{k,p}^f\},\; \text {where}\;\; \boldsymbol x_{k,0}^f := \boldsymbol x_k\\
 \label {sample:2}\mathcal{B}_{\epsilon,k}^{h^i} &:= \{\boldsymbol x_{k,0}^{h^i},\boldsymbol x_{k,1}^{h^i},\cdots,\boldsymbol x_{k,p}^{h^i}\},\; \text{where}\;\; \boldsymbol x_{k,0}^{h^i} := \boldsymbol x_k\\
 \label {sample:3}\mathcal{B}_{\epsilon,k}^{g^j} &:= \{\boldsymbol x_{k,0}^{g^j},\boldsymbol x_{k,1}^{g^j},\cdots,\boldsymbol x_{k,p}^{g^j}\},\;\text{where}\;\; \boldsymbol x_{k,0}^{g^j} := \boldsymbol x_k
\end{align}
are sets of $p$ (sample size) independent and identically distributed random points uniformly sampled from
\begin{align}
\label {sample:4}\mathbb{B}_\epsilon(\boldsymbol x_k) :={\lVert \boldsymbol  x-\boldsymbol x_k\rVert_2 \leq \epsilon},
\end{align}
where $\epsilon$ is the sample radius.
To indicate the progress of the algorithm iterations, an infeasibility vector is defined as
\begin{equation} \label{infeas}
\begin{aligned}
\hspace*{-1.2cm}
\boldsymbol \sigma_k(\boldsymbol x_k):= & \left(\begin{array}{c}
     \boldsymbol |\boldsymbol h(\boldsymbol x_k)|  \\
      \max(\boldsymbol g(\boldsymbol x_k)-\overline{\boldsymbol g},0)\\  
      \max\big(\underline{\boldsymbol g}-\boldsymbol g( \boldsymbol x_k),0)
\end{array} 
\right).
\end{aligned}
\end{equation}
The following model reduction is also defined in terms of primal and dual infeasibility, which can be zero only if $\boldsymbol  x_k$ is $\epsilon$-stationary\cite{curtis2012sequential}:
\begin{align} \label{reduction}
\hspace*{-0.3cm}\Delta \boldsymbol q_k :&= \rho f( \boldsymbol x_k)+\boldsymbol\sum{ \boldsymbol \sigma_k( \boldsymbol x_k)}-\rho \max_{x \in \mathcal{B}_{\epsilon,k}^{f}}\{f(\boldsymbol x_k)+\nabla f(\boldsymbol x)^\top\boldsymbol d_k\} \notag \\
-&\frac{1}{2}\boldsymbol d_k^\top\boldsymbol H_k\boldsymbol d_k -\sum_{i=1}^{m_h}\max_{x \in \mathcal{B}_{\epsilon,k}^{h^i}}|h^i(\boldsymbol x_k)+\nabla h^i(\boldsymbol x)^\top\boldsymbol d_k| \notag \\
-&\sum_{j=1}^{m_g}\max_{x \in \mathcal{B}_{\epsilon,k}^{g^j}}\big\{\max\{g^j(\boldsymbol x_k)+\nabla g^j(\boldsymbol x)^\top\boldsymbol d_k-\overline{g}^j,0\}\big\} \notag \\
-&\sum_{j=1}^{m_g}\max_{x \in \mathcal{B}_{\epsilon,k}^{g^j}}\big\{\max\{-g^j(\boldsymbol x_k)-\nabla g^j(\boldsymbol x)^\top\boldsymbol d_k+\underline{g}^j,0\}\big\}.
\end{align}
Then the SQP-GS algorithm is presented in Algorithm \ref{algo1}.

The SQP-GS algorithm generalizes the traditional SQP method to nonsmooth constrained problems. 
When solving a smooth constrained problem, the sample size $p$ can be chosen as zero, in which case the SQP-GS algorithm reduces to the traditional SQP method. 

If a function is known to be continuously differentiable everywhere in $\mathbb{R}^n$, it is unnecessary to sample its gradient at nearby points. 
This can significantly improve the performance of the algorithm because the evaluation of  the linear inequality constraints from the quadratic programming can be largely eliminated. 
Besides, for those functions that depend on fewer than $n$ variables, sampling fewer points can still yield good results. 
Also, since the points are sampled independently, it allows to use parallel computing  to further reduce CPU time.

\begin{algorithm}[!ht]
\caption{SQP-GS Algorithm}\label{algo1}
\begin{algorithmic}[1]
\State \textbf{Initialization}: Set $k=1$, $K_{\max}=100$; Choose a  sampling radius $\epsilon>0$, penalty parameter $\rho>0$,  constraint violation tolerance $\tau>0$,  sample size $p > 0$, line search constant $\varpi \in (0,1)$, backtracking constant $\gamma\in (0,1)$, sampling radius reduction factor $\mu_\epsilon\in(0,1) $, penalty parameter reduction factor $\mu_\rho  \in (0,1)$, constraint violation tolerance reduction factor $\mu_\tau \in (0,1)$, infeasibility tolerance  $\nu_{in}>0$, and stationarity tolerance parameter $\nu_s>0$. Choose an initial iterate $\boldsymbol x$.
\State \textbf{while} $k<K_{\max}$
\State \textbf{if} $\max(  \Delta \boldsymbol q_k)<\nu_s$ and $\max(\boldsymbol \sigma_k)<\nu_{in}$\\
  \hspace {0.5cm}\ Output solution and stop.
\State \textbf{end if}
\State \textbf{Gradient sampling}: Generate $\mathcal{B}_{\epsilon,k}^f$, $\mathcal{B}_{\epsilon,k}^h$, and $\mathcal{B}_{\epsilon,k}^g$ by (\ref{sample:1})--(\ref{sample:4}).
\State \textbf{Search direction calculation}: Solve (\ref{model:1})--(\ref{model:7}) to get $(\boldsymbol d_k,z_k,\overline{\boldsymbol r}_k ,\underline{\boldsymbol r}_k,\boldsymbol e_k)$.
\State \textbf{L-BFGS update}: Limited-memory Broyden-Fletcher-Goldfarb-Shanno (L-BFGS) method\cite{nocedal1980updating} that is typical in smooth optimization is used to update $\boldsymbol H_k$.
\State  \textbf{if} $\Delta \boldsymbol q_k>\nu_s \epsilon^2$\\
           \hspace {0.5cm}\textbf{go} to step 20.\\
       \textbf{else} \\
           \hspace {0.5cm}\textbf{if}  $\max(\boldsymbol \sigma_k) \leq \tau$ \\
           \hspace {0.5cm}\hspace {0.5cm} \textbf{set} $\tau\leftarrow \mu_\tau \tau$\\
           \hspace {0.4cm} \textbf{else}\\
               \hspace {0.5cm}\hspace {0.5cm}\textbf{set} $\rho \leftarrow \mu_\rho \rho$\\
           \hspace {0.5cm}\textbf{end if}\\
           \hspace {0.4cm} \textbf{set} $\epsilon \leftarrow \mu_\epsilon\epsilon$, $\beta_k\leftarrow 0$\\
           \hspace {0.4cm} \textbf{go} to step 21\\
            \textbf{end if}
\State \textbf{Line search}: Set $\beta_k$ as the largest value in the sequence $\{1,\gamma,\gamma^2,\cdots\}$ such that $\boldsymbol x_{k+1} \leftarrow \boldsymbol x_k+ \beta_k \boldsymbol d_k$  satisfies:
  \begin{equation}
 \begin{aligned}
\hspace*{-0.5cm}\rho f(\boldsymbol x_{k+1})+\sum{\boldsymbol \sigma_k(\boldsymbol x_{k+1}}) \leq &\rho f(\boldsymbol x_k)+\sum{\boldsymbol \sigma_k(\boldsymbol x_k)}\\
           &-\varpi \beta_k\Delta \boldsymbol q_k
 \end{aligned}
\end{equation}
  \State \textbf{Iteration increment}: Set $k \leftarrow k+1$
\State \textbf{end do}
	\end{algorithmic}
\end{algorithm}


\section{Solving SSSC-OPF by SQP-GS}  \label{sss-opf}

Here we discuss how to apply the SQP-GS algorithm in Section \ref{sqp-gs} to solve the SSSC-OPF problem.

\subsection{Model of SSSC-OPF}

The model of SSSC-OPF is actually a `standard' OPF model  defined as a smooth nonlinear programming problem with an additional small-signal stability constraint. 
Specifically, the SSSC-OPF model can be represented as follows.
\begin{enumerate}
   \item  Minimizing the generation cost is usually considered as the objective function
    \begin{equation}
      \label {SSS:1}   f = \sum_{i \in S_G}(a_{2i}P_{\textrm{G}i}^2+a_{1i}P_{\textrm{G}i}+a_{0i} ),\quad 
    \end{equation}
   where  $P_{\textrm{G}i}$ is the active power output of the $i$th generator, and $a_{2i}$, $a_{1i}$, and $a_{0i}$ are the corresponding cost coefficients. 
  \item Power flow equations for bus $i\in S_\textrm{B}$
        \begin{eqnarray}
        \label{power_equ}
        P_{\textrm{G}i}-P_{\textrm{L}i}-\sum_{j \in i}V_i Y_{ij}V_j \cos(\theta_i-\theta_j-\alpha_{ij} ) = 0 \\
       Q_{\textrm{G}i}-Q_{\textrm{L}i}-\sum_{j \in i}V_i Y_{ij}V_j \sin(\theta_i-\theta_j-\alpha_{ij} ) = 0,
        \end{eqnarray}
        where  $Q_{\textrm{G}i}$ is the reactive power output of the $i$th generator. 
   \item Initial condition equations for generator $i \in S_\textrm{G}$ with a two-axis model.
   \begin{itemize}
  \item Stator algebraic equations
      \begin{eqnarray}
   \hspace*{-0.8 cm}   E_{\textrm{d}i}^{'}-V_i\sin(\delta_i-\theta_i)-R_{\textrm{s}i}I_{\textrm{d}i}+X_{\textrm{q}i}^{'}I_{\textrm{q}i} = 0 \\
    \hspace*{-0.8cm}  E_{\textrm{q}i}^{'}-V_i\cos(\delta_i-\theta_i)-R_{\textrm{s}i}I_{\textrm{q}i}-X_{\textrm{d}i}^{'}I_{\textrm{d}i} = 0.
      \end{eqnarray}
      \item The generator terminal power can be obtained by the product of voltage and current transformed from $\textrm{d}-\textrm{q}$ rotor reference frame to network reference frame as
      \begin{align}
   \hspace*{-0.8cm} P_{\textrm{G}i}-V_iI_{\textrm{d}i}\sin(\delta_i-\theta_i)-V_iI_{\textrm{q}i}\cos(\delta_i-\theta_i) &= 0 \\
   \hspace*{-0.9cm}\, Q_{\textrm{G}i}-V_iI_{\textrm{d}i}\cos(\delta_i-\theta_i)+V_iI_{\textrm{q}i}\sin(\delta_i-\theta_i) &= 0.
  \end{align}
    \item In steady state, ${dE_{\textrm{q}i}^{'}}/{dt} $ in (\ref{gen:Eq}) and ${dE_{\textrm{d}i}^{'}}/{dt} $ in (\ref{gen:Ed})   are zero and thus
  \begin{eqnarray}
     E_{\textrm{fd}i}-E_{\textrm{q}i}^{'}-(X_{\textrm{d}i}-X_{\textrm{d}i}^{'})I_{\textrm{d}i} = 0\\
       E_{\textrm{d}i}^{'}-(X_{\textrm{q}i}-X_{\textrm{q}i}^{'})I_{\textrm{q}i} = 0.
  \end{eqnarray}
 \end{itemize}
    \item Technical constraints include
\begin{align}
\underline{V}_i \leq V_i &\leq \overline{V}_i  &i \in S_\textrm{B} \\
\underline{P}_{\textrm{G}i} \leq P_{\textrm{G}i} &\leq \overline{P}_{\textrm{G}i}   &i \in S_\textrm{G} \\
\underline{Q}_{\textrm{G}i} \leq Q_{\textrm{G}i} &\leq \overline{Q}_{\textrm{G}i}  &i \in S_\textrm{G} \\
  \label {SSS:13}  I_{ij}^2 &\leq \overline{I}_{ij}^2  &(i,j) \in S_\textrm{Line},
\end{align}
 where $I_{ij}$ is the current of line $(i,j)$, $S_\textrm{Line}$ is the set of all lines and $\overline{(\cdot)}$ and $\underline{(\cdot)}$ denote the upper and lower limits.

\item Small-signal stability constraint:
\begin{align}
   \label {SSS:14}   \eta(\boldsymbol  x) \leq \overline{\eta}
\end{align}
\end{enumerate}
where $ \boldsymbol x =\big[\boldsymbol P_\textrm{G}^\top\;\boldsymbol Q_\textrm{G}^\top\;\boldsymbol V^\top\;\boldsymbol \theta^\top\;\boldsymbol \delta^\top\;\boldsymbol E_\textrm{d}^{'\top}\;\boldsymbol E_\textrm{q}^{'\top}\;\boldsymbol I_\textrm{d}^\top\;\boldsymbol I_\textrm{q}^\top\;\boldsymbol E_{\textrm{fd}}^\top\big]^\top$ is the vector of variables in the model. The choice of $\overline{\eta}$ depends on the system characteristics and can be determined based on offline stability studies.

\subsection{Employing SQP-GS to Solve SSSC-OPF}

 Obviously, the SSSC-OPF model belongs to the type of optimization problem in (\ref{model:nonlin}). 
 Before employing the SQP-GS method that relies on the gradients to construct the QP subproblem, 
 the gradients of all functions in the model with respect to $\boldsymbol x$, the variables of the model, should be derived. 
 The objective function and constraint functions in (\ref{power_equ})--(\ref{SSS:13}) are smoothly nonlinear or linear, 
 and their gradients with respect to $ \boldsymbol x$ can be easily derived. 
 As discussed in Section \ref{sss}, the  spectral abscissa function in (\ref{SSS:14})  is implicit but the function value and its gradient can also be evaluated. 
 
 As discussed in Section \ref{sqp-gs}, it is unnecessary to sample the gradients of the smooth functions. 
 In SSSC-OPF, the function in (\ref{SSS:14}) is the only nonsmooth function, and thus only its gradient need to be sampled, 
 which requires the following steps for the $k$ iteration: 
\begin{enumerate}
\item \textbf{Sampling Points}: Generate $p$ points $\mathcal{B}_{\epsilon,k}^{\eta}$ by (\ref{sample:4}).
\item \textbf{Eigenvalue Analysis}: Set up state matrix  $\boldsymbol A$ for each point $\boldsymbol x \in  \mathcal{B}_{\epsilon,k}^{\eta} $ and calculate the most critical eigenvalue $\lambda_\eta$ and the corresponding left and right eigenvectors $\boldsymbol \psi_\eta$ and $\boldsymbol \phi_\eta$ for each $\boldsymbol A$;  
\item  \textbf{Sensitivities}: Calculate the gradients ${\partial \tilde{\boldsymbol A }}/{\partial  x_i}$, ${\partial \tilde{\boldsymbol B }}/{\partial  x_i}$, ${\partial \tilde{\boldsymbol C }}/{\partial  x_i}$, and ${\partial \tilde{\boldsymbol D }}/{\partial  x_i}$  for each point with respect to the $i$th variable $x_i \in \boldsymbol x $;
obtain the gradients  ${\partial\boldsymbol A }/{\partial  x_i}$ for each $\boldsymbol A$ in terms of (\ref{eigen:A}); 
and calculate the eigenvalue sensitivity ${\partial\lambda_\eta}/{\partial x_i}$ in (\ref{eq:eigen_sen}) for  $\boldsymbol x \in  \mathcal{B}_{\epsilon,k}^{\eta} $. 
\item  \textbf{Gradient}: Get the gradient ${\partial \eta}/{\partial x_i}$  with respect to the $i$th variable $x_i \in \boldsymbol x $  referred to (\ref{eigen:sa})  for each point.
\end{enumerate}
Note that the steps (2)--(4) can be performed in parallel.

\section{Case Studies} \label{result}

The proposed method is applied to the  WSCC 3-machine 9-bus , New England  10-machine 39-bus, and the modified IEEE 57-machine 118-bus systems 
to illustrate the effectiveness in solving SSSC-OPF. For all systems, the generators are described by the two-axis model with an IEEE type-I exciter. 
The loads are modeled as constant power. 

The SQP-GS method is implemented in MATLAB by using CPLEX 12.60 \cite{cplex} as the QP solver for the subproblem. 
The eigenvalues and eigenvectors are computed by QR decomposition using the MATLAB function \emph{eig}. 
Flat start is used for which all voltage angles are set to be zero, all voltage magnitudes are set to be $1.0$ p.u., 
$P_G = (\overline{P}_G+\underline{P}_G)/2$, and $Q_G = (\overline{Q}_G+\underline{Q}_G)/2$. 
The parameters of SQP-GS are chosen as $\rho= 0.1$, $\mu_\rho= 0.5$, $\epsilon= 0.1$, $\mu_\epsilon = 0.5$, $\tau= 0.1$, $\mu_\tau=0.8$, $\varpi = 1$, $\gamma= 0.8$ from \cite{curtis2012sequential}. 
The tolerances $\nu_{in}$ and $\nu_{s}$ are set to be $10^{-3}$ and $10^{-2}$, respectively.

\subsection{WSCC 3-Machine 9-Bus System}

The WSCC 9-bus system is often used for stability analysis \cite{sauer1998power}. The system data, including power limit and security data can be found  \cite{6531914}. 
The generator cost coefficients are listed in Table \ref{table_1}. To analyze the effectiveness of the proposed method, we consider the following three cases:
\begin{enumerate}
\item \textit{Case 0}: Base case, which is a standard OPF without any small-signal stability constraint.
\item \textit{Case 1}: SSSC-OPF without any binding small-signal stability constraint.
\item \textit{Case 2}: SSSC-OPF with a binding small-signal stability constraint.
\end{enumerate}
\begin{table}[htb]
\caption{Generation Cost Coefficients for WSCC 9-Bus System}
\label{table_1}
\begin{center}
{
\begin{tabular}{c c c c }\hline
Generator\# &$a_2$,\,\$/$(\text{MW})^2$ &$a_1$,\,\$/$\text{MW}$   &$a_0$,\,\$ \\\hline
1     & $9.76\times 10^{-4}$       & 14.712      & 0.00\\
2     & $7.20\times 10^{-4}$       & 11.290      & 0.00\\
3     & $5.46\times 10^{-4}$       & 8.001       & 0.00\\\hline
\end{tabular}
}
\end{center}
\end{table}

\begin{itemize}
\item \textbf{Case 0}: We use IPM to solve the standard OPF without small-signal stability constraint
and the results are listed in the first row of Table \ref{table_2}. 
The eigenvalue analysis is performed and the spectral abscissa $\eta$ is $-0.04$.

\item \textbf{Case 1}: A security margin $\overline \eta = -0.01$ is used, which is larger than the spectral abscissa in the base case. 
Since the small-signal stability constraint is not binding in this case, the problem can be successfully solved by the SQP method without gradient sampling and 
the results are the same as those in the base case. 

\item \textbf{Case 2}: $\overline \eta$ is set to be three different values, all of which are less than the spectral abscissa in the base case. 
When $\overline \eta$ is set to be $-0.45$ or $-0.5$, the SQP will not converge without a GS procedure. 
Here the sample size is chosen as $10$. 
From Table \ref{table_2} it is seen that the power outputs of the generators are re-dispatched to satisfy the small-signal stability constraint. 
Also, the more binding the small-signal stability constraint, the more expensive the generation cost is. 
As $\overline{\eta}$ decreases, the active power generated by generator G3 which has the cheapest generation cost gradually goes down, 
mainly because G3 has the smallest inertial constant and generating more power from the other two generators can help improve the stability. 
When $\overline \eta=-0.5$, the re-dispatch from SSSC-OPF will allow significantly more generation from G1 which has the most expensive generation cost but the largest inertial constant.
\end{itemize}

\begin{table*}[htb]
\caption{Summary of SSSC-OPF Results With Different Small-Signal Constraints for WSCC 3-Machine 9-Bus System}
\label{table_2}
\begin{center}
{
\begin{tabular}{c c c c c c c c c c}\hline
         &Generation Cost ($\$/h$)    &$\eta$   &$P_{\textrm{G}1}(\text{MW})$  &$P_{\textrm{G}2} (\text{MW})$ &$P_{\textrm{G}3} (\text{MW})$ &$V_1 (\text{pu})$ &$V_{2} (\text{pu})$  &$V_{3} (\text{pu})$\\\hline
Base Case           &2901.3       &-0.04       &25.0   &25.0          &276.0      &1.040  &1.045     &1.022 \\ 
Case 1 ($\eta \leq -0.01$)   &2901.3       &-0.04       &25.0   &25.0          &276.0      &1.040  &1.045     &1.022\\
Case 2 ($\eta \leq -0.30$)   &2919.7       &-0.30       &25.0   &47.1          &252.4      &1.040  &1.045     &1.045\\
Case 2 ($\eta \leq -0.45$)   &2957.5       &-0.45       &25.0   &59.8          &239.3      &1.026  &1.045     &1.045\\
Case 2 ($\eta \leq -0.50$)   &3217.9       &-0.50       &84.4   &25.0          &211.8      &0.997  &1.045     &1.022\\
\hline
\end{tabular}
}
\end{center}
\end{table*}

\subsection{New England 10-Machine 39-Bus System}

The full dynamic data of the New England 10-machine 39-bus system are extracted from \cite{yeu2010small} and the economic and technical data are from \cite{5593195}. 
For the voltage magnitude limits, we choose $\overline{V}= 1.1$ p.u. and $\underline{V}=  0.9$ p.u. for all buses. 
For standard OPF, the spectral abscissa $\eta = -0.11$.
To reduce the spectral abscissa, the SSSC-OPF is applied with $\overline{\eta} = -0.2$.

\begin{figure}
    \centering
    \includegraphics[width = 0.48\textwidth]{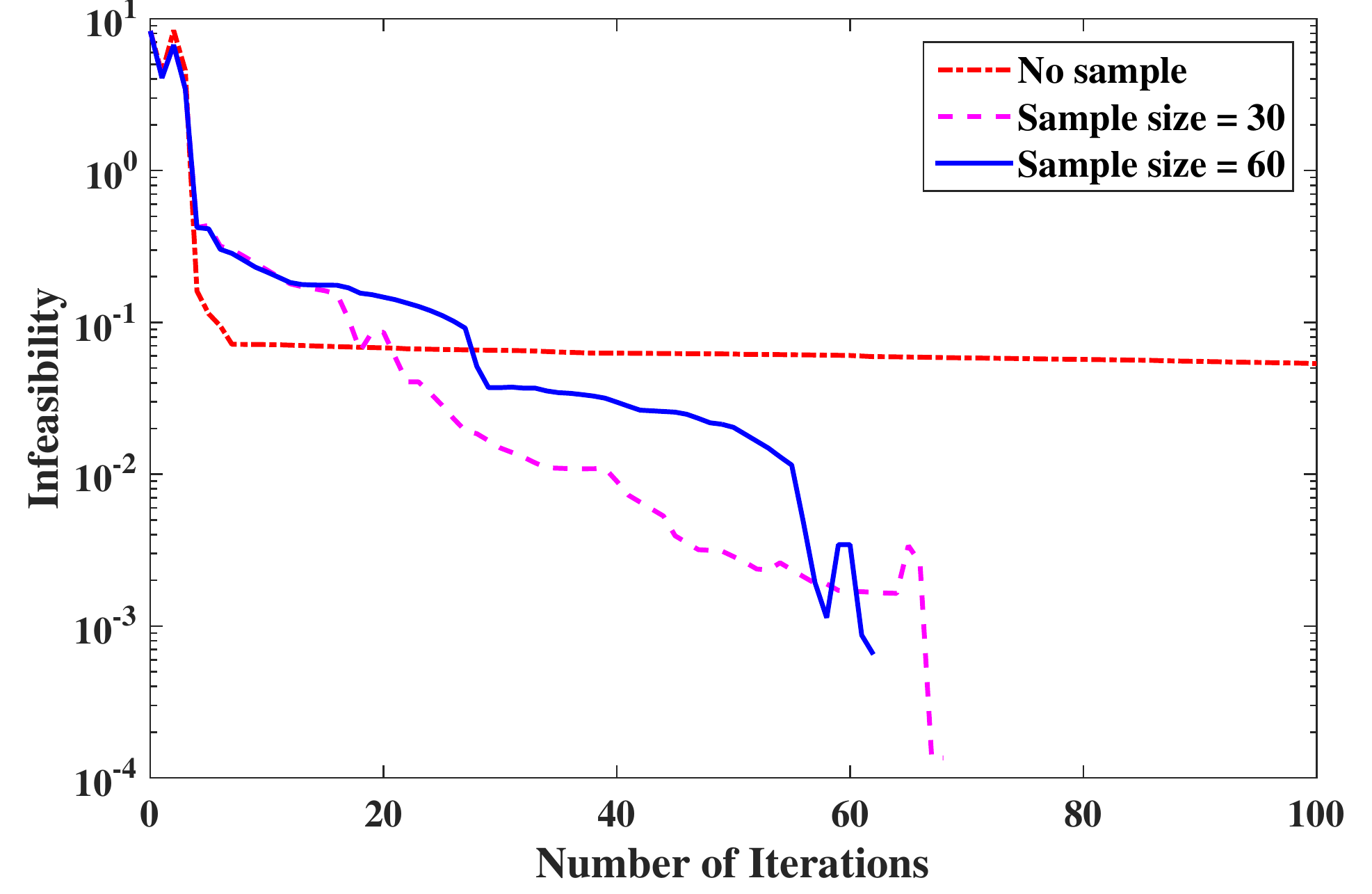}
    \caption{Infeasibility curve for IEEE 39-bus system.}
    \label{fig:IEEE39_infeasibility}
\end{figure}

\begin{figure}
    \centering
    \includegraphics[width = 0.48\textwidth]{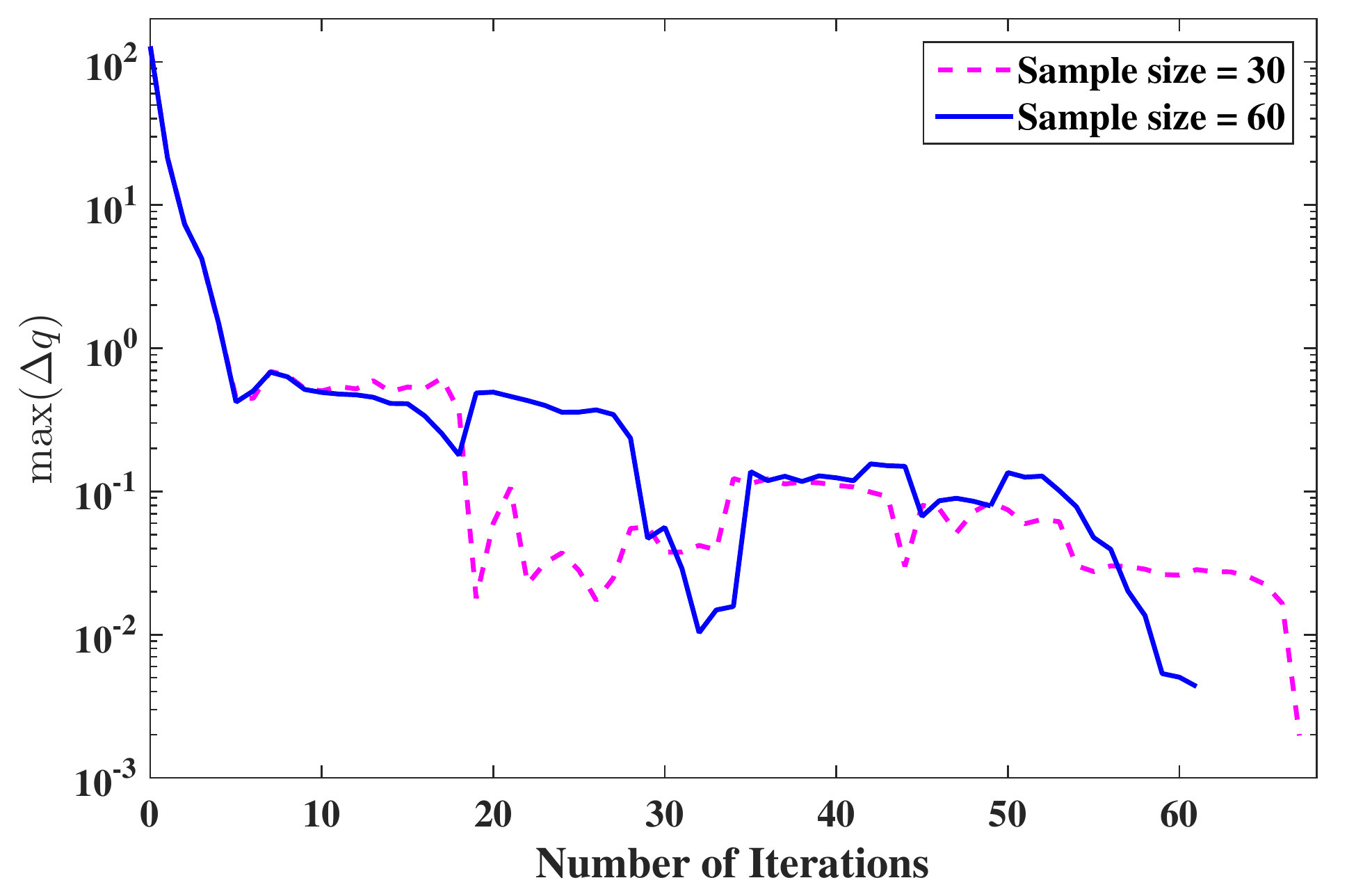}
    \caption{Model reduction curve for IEEE 39-bus system.}
    \label{fig:modify_reduction}
\end{figure}

\begin{figure}
    \centering
    \includegraphics[width = 0.48\textwidth]{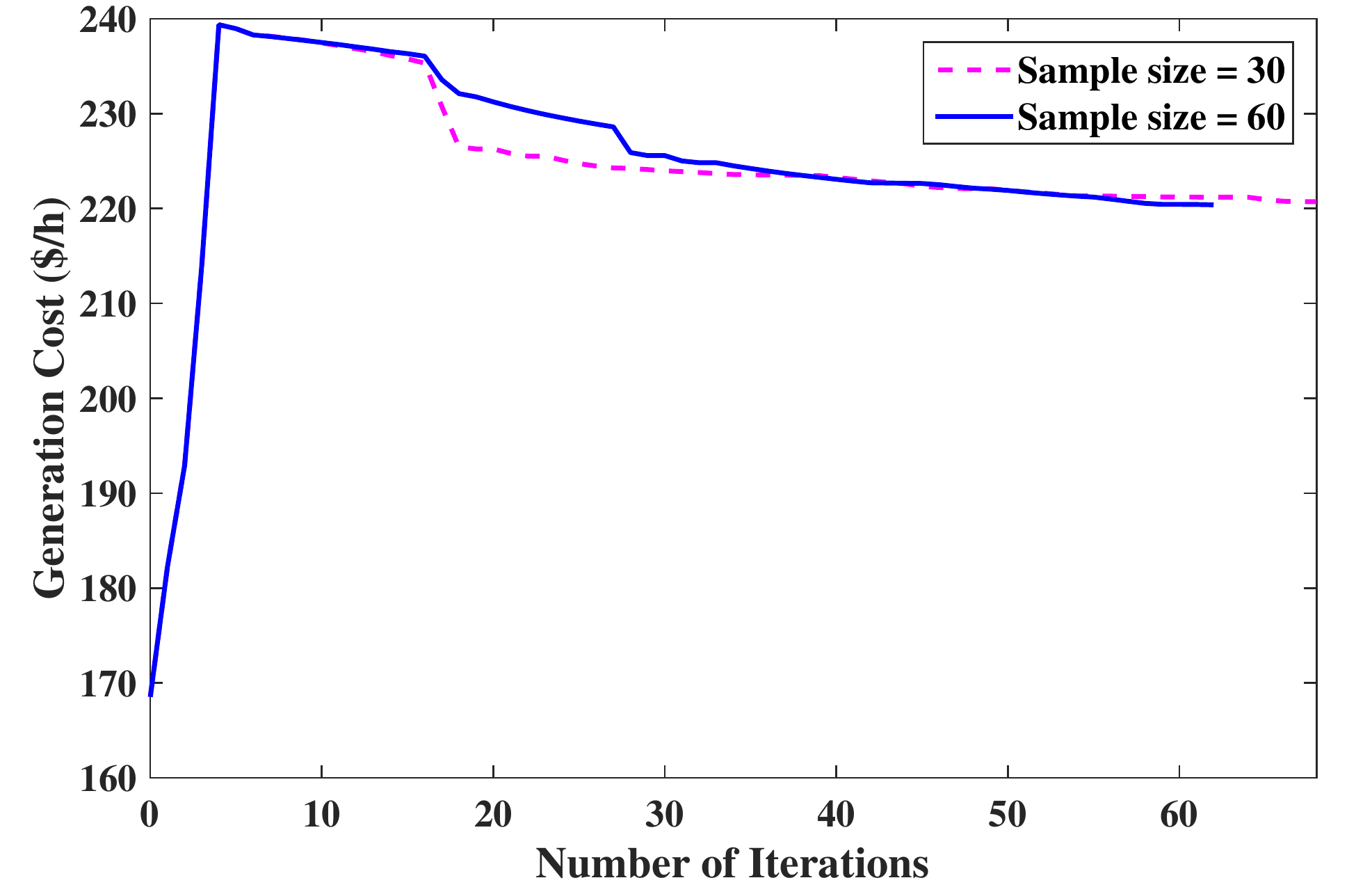}
    \caption{Objective curve for IEEE 39-bus system.}
    \label{fig:obj}
\end{figure}

The SSSC-OPF is calculated by SQP-GS with no samples, $30$ samples, and $60$ samples, respectively. 
As in Fig. \ref{fig:IEEE39_infeasibility}, the infeasibility defined in (\ref{infeas}) cannot reduce to the tolerance with $100$ iterations when the sample size is $0$. 
By contrast, when the sample size is $30$ or $60$, the infeasibility reduces to $v_{in}=10^{-3}$ rapidly. 
Also, as shown in Fig. \ref{fig:modify_reduction}, the model reduction defined in (\ref{reduction}) decreases to an acceptable tolerance with sampling gradients. 
In Fig. \ref{fig:obj} we show the generation cost which stably approaches an optimal value. 
From Figs. \ref{fig:IEEE39_infeasibility}--\ref{fig:obj}, we can see that the gradient sampling plays an important role in solving the SSSC-OPF problems. 
Moreover, based on our tests, the GS with 30 samples are good enough to improve the search direction of SQP for the SSSC-OPF problem.

We also compare SQP-GS with the numerical eigenvalue sensitivity based IPM (IPM-NES) \cite{5593195} 
and the results are listed in Table \ref{table_cost_compare}, where $\overline{P}_\textrm{G}$ is the maximum active power output of a generator. 
The results for the OPF without small signal stability constraint (denoted by `OPF') is also listed for reference. 
IPM-NES can also solve the problem but needs higher cost than SQP-GS to ensure the same level of small-signal stability. 
Comparing the generations of  IPM-NES and SQP-GS, we can see that the generations of G3 and G5 are the same while
there are significant differences for the other generators. 
Actually the active power outputs of G3 and G5 always reach the maximum power output in all methods, which is mainly due to their much cheaper generation cost. 

The critical modes that change in the first $10$ iterations for both methods are listed in Table \ref{table_crical_mode}. 
We can see that after $5$ iterations SQP-GS moves to the binding critical modes for small-signal stability constraint and stays in this mode until convergence. 
However, IPM-NES gets the same mode at the $3$rd and $10$th iteration but suffers from oscillations between some critical modes during iterations. 
\begin{table}[H]
\caption{Comparison Between SQP-GS and IPM-NES for 39-Bus System}
\label{table_cost_compare}
\begin{center}
\renewcommand{\arraystretch}{1.3}
{
\begin{tabular}{c c c c c}\hline
   & \tabincell{c}{SQP-GS \\ (MW)}       & \tabincell{c}{IPM-NES \\ (MW)}     & \tabincell{c}{OPF \\ (MW)}   & \tabincell{c}{$\overline{P}_\textrm{G}$ \\ (MW)} \\
   \hline
G1  &0.00        &402.5  &0.00    &402.5\\
G2  &698.8       &132.4  &747.5    &747.5\\
G3  &920.0       &920.0  &920.0   &920.0\\
G4  &298.8       &100.8  &0.00    &862.5\\
G5  &747.5       &747.5  &747.5   &747.5\\
G6  &740.0       &479.8  &862.5  &862.5\\
G7  &342.3       &862.5  &0.00   &862.5\\
G8  &557.1       &181.8  &805.0  &805.0\\
G9  &804.3       &1000.4 &883.3  &1035.0\\
G10 &1047.1      &1308.8 &1195.0 &1380.0\\\hline
Cost\,(\$/h)&220.7     &239.1 &213.9 &--\\\hline
\end{tabular}
}
\end{center}
\end{table}
\subsection{IEEE 54-Machine 118-Bus System}

We also test the proposed method on a modified version of the IEEE 54-machine 118-bus benchmark system with dynamic data from \cite{IEEE_118}.  
This system has $54$ synchronous machines with IEEE type-$1$ exciters, $20$ of which are synchronous compensators used only for reactive power support and $15$ of which are motors. 
For a standard OPF, the spectral abscissa $\eta =0.35$. We set $\overline{\eta}$ to be $-0.1$. The voltage magnitudes at all buses must be between $0.9$ p.u. and $1.1$ p.u. . 
The reactive power limits can be found in \cite{IEEE_118_limits}.

IPM-NES fails to solve this problem. By contrast, SQP-GS only needs $13$ and $12$ iterations to solve the problem, respectively, when there are $30$ samples and $60$ samples. 
The generation cost of SSSC-OPF will increase to $6110.8\$/$h  from $5779.2\$/$h for the standard OPF with $\eta =0.35$. 

\begin{figure}
    \centering
    \includegraphics[width = 0.48\textwidth]{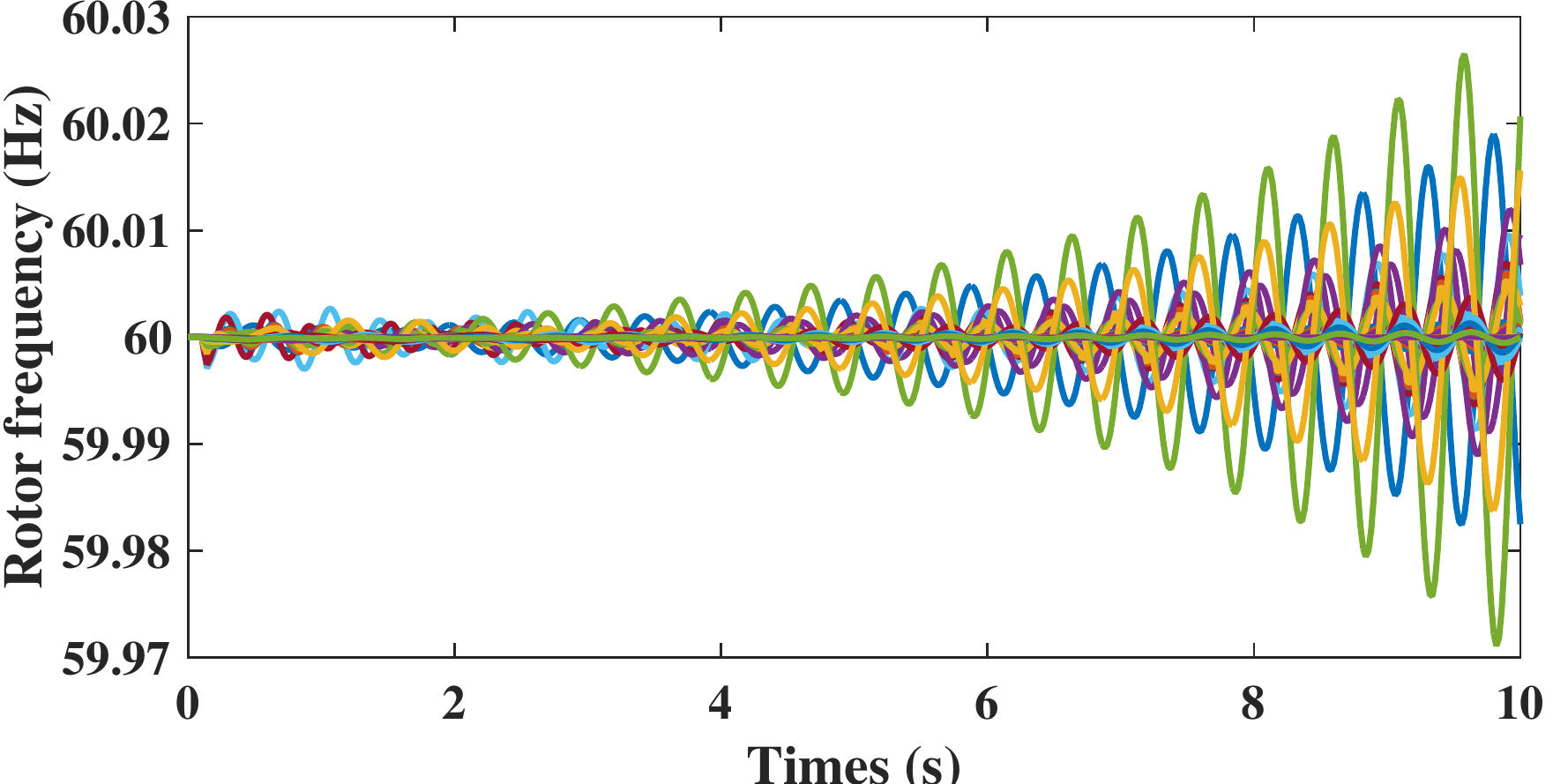}
    \caption{Time-domain simulation of 118-bus system for standard OPF.}
    \label{fig:unstable}
\end{figure}

\begin{figure}
    \centering
    \includegraphics[width = 0.48\textwidth]{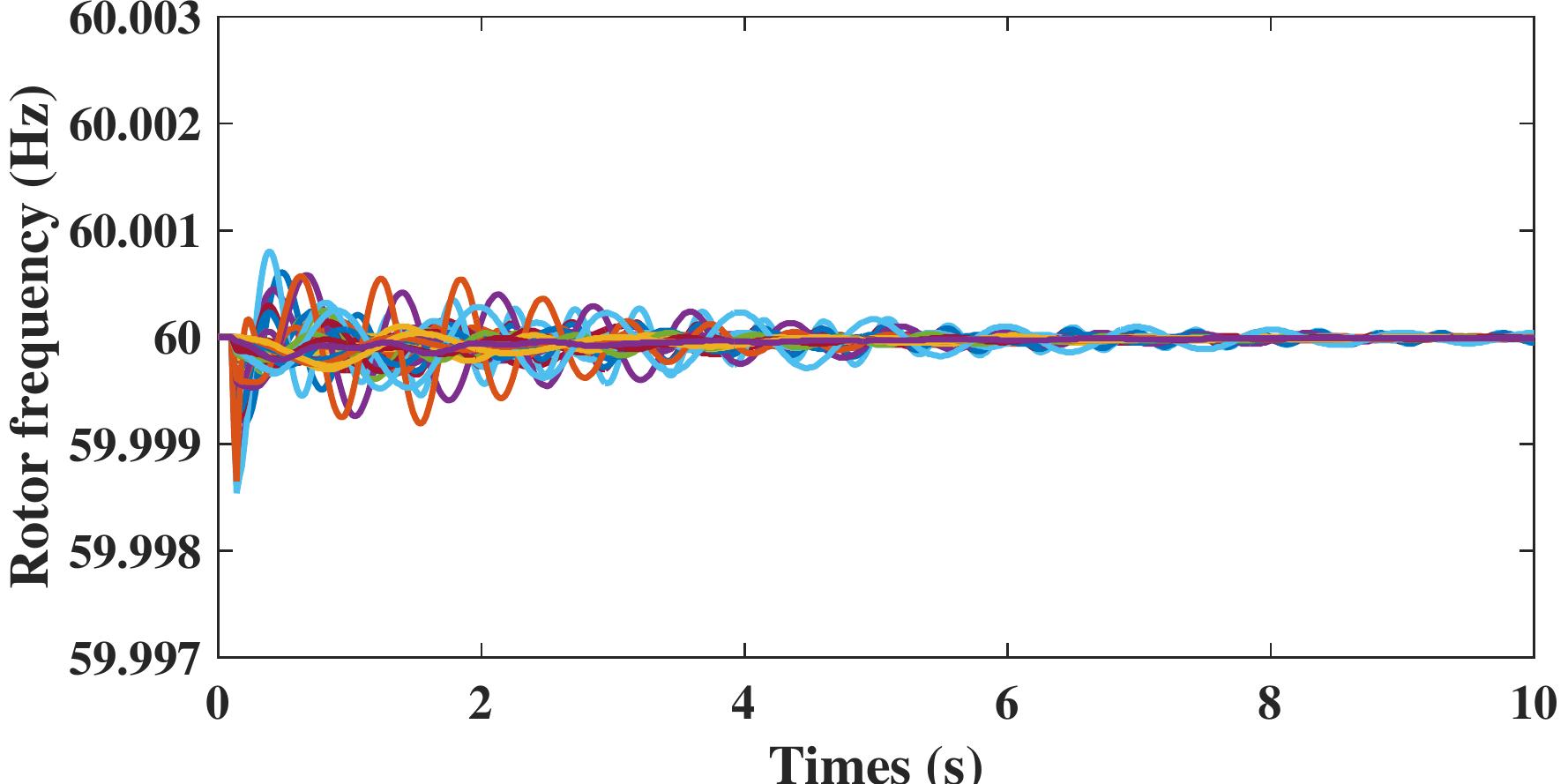}
    \caption{Time-domain simulation of 118-bus system for SSSC-OPF.}
    \label{fig:stable}
\end{figure}

Since eigenvalue analysis and the SSSC-OPF use the linearized system model, they cannot take into account the full nonlinear behavior of the power system. 
Therefore, in order to further validate the solution of SSSC-OPF from SQP-GS, we increase the load at bus $2$ by $5$ MW and perform
time-domain simulation using the full nonlinear power system model for the operating states obtained from both the standard OPF and SSSC-OPF.
The rotor frequencies of all generators for the standard OPF and SSSC-OPF are shown in Figs. \ref{fig:unstable} and \ref{fig:stable}, respectively. It is seen that under the standard OPF the system is unstable while under SSSC-OPF the oscillation can be quickly damped. 
\begin{table}[H]
\caption{Critical Modes in 10 Iterations for 39-Bus System}
\label{table_crical_mode}
\begin{center}
{
\begin{tabular}{c c c }\hline
Iter. &SQP-GS    &IPM-NES \\\hline
1         &\tabincell{l}{\textbf{-0.1193 $\pm$ j0.4060}\\ -0.1389 $\pm$ j0.3919\\  -0.1664 $\pm$ j0.4264\\  -0.1773 $\pm$ j4.1030} 
          &\tabincell{l}{\textbf{-0.0965 $\pm$ j 0.5086}\\ -0.1537 $\pm$ j0.5302\\  -0.1869 $\pm$ j0.7166\\ -0.1876 $\pm$ j6.2010}\\\hline
2         &\tabincell{l}{\textbf{-0.1394 $\pm$ j0.4056}\\ -0.1597 $\pm$ j0.4263\\  -0.1706 $\pm$ j0.4484\\ -0.1836 $\pm$ j4.3023}    
          &\tabincell{l}{\textbf{-0.1860 $\pm$ j0.7117}\\ -0.1935 $\pm$ j3.9960\\  -0.1965 $\pm$ j7.5179 \\-0.2010 $\pm$ j6.0230} \\\hline
3         &\tabincell{l}{\textbf{-0.1578 $\pm$ j0.4271}\\ -0.1799 $\pm$ j0.4217 \\ -0.1865 $\pm$ j4.4437\\ -0.1994 $\pm$ j0.4592}
          &\tabincell{l}{\textbf{-0.1776 $\pm$ j6.1767}\\ -0.1888 $\pm$ j3.9706 \\ -0.2102 $\pm$ j7.5391\\ -0.2104 $\pm$ j5.4174}\\\hline
4         &\tabincell{l}{\textbf{-0.1762 $\pm$ j0.4374}\\ -0.1837 $\pm$ j4.5249 \\ -0.1842 $\pm$ j0.4306\\ -0.1987 $\pm$ j6.9051}
          &\tabincell{l}{\textbf{-0.1768 $\pm$ j4.0281}\\ -0.1917 $\pm$ j6.1061 \\ -0.2035 $\pm$ j7.5423 \\-0.2069 $\pm$ j0.6721}\\\hline  
5         &\tabincell{l}{\textbf{-0.1872 $\pm$ j6.8941}\\ -0.1895 $\pm$ j0.4458 \\ -0.1904 $\pm$ j4.4475\\ -0.1970 $\pm$ j0.4421}
          &\tabincell{l}{\textbf{-0.1815 $\pm$ j4.0092}\\ -0.1945 $\pm$ j6.1100 \\ -0.2036 $\pm$ j7.5421 \\-0.2066 $\pm$ j0.6721}\\\hline 
6         &\tabincell{l}{\textbf{-0.1829 $\pm$ j6.6878}\\ -0.1921 $\pm$ j4.2950 \\ -0.1948 $\pm$ j0.4488\\ -0.2021 $\pm$ j0.4487}
          &\tabincell{l}{\textbf{-0.1865 $\pm$ j3.9870}\\ -0.1918 $\pm$ j6.1391 \\ -0.2031 $\pm$ j7.5365 \\-0.2069 $\pm$ j0.6708}\\\hline 
7         &\tabincell{l}{\textbf{-0.1835 $\pm$ j6.5399}\\ -0.1932 $\pm$ j4.2018 \\ -0.1938 $\pm$ j0.4522\\ -0.2053 $\pm$ j0.4516}
          &\tabincell{l}{\textbf{-0.1873 $\pm$ j6.1417}\\ -0.1922 $\pm$ j3.9662 \\ -0.2032 $\pm$ j7.5364 \\-0.2066 $\pm$ j0.6708}\\\hline  
8         &\tabincell{l}{\textbf{-0.1835 $\pm$ j6.5089}\\ -0.1938 $\pm$ j4.1797\\  -0.1941 $\pm$ j0.4527\\  -0.2060 $\pm$ j0.4525}
          &\tabincell{l}{\textbf{-0.1803 $\pm$ j4.0352}\\ -0.1887 $\pm$ j0.7124 \\ -0.1974 $\pm$ j7.5183 \\-0.2015 $\pm$ j6.0104}\\\hline  
9         &\tabincell{l}{\textbf{-0.1838 $\pm$ j6.4765}\\ -0.1943 $\pm$ j0.4533\\  -0.1946 $\pm$ j4.1570\\ -0.2067 $\pm$ j0.4534}
          &\tabincell{l}{\textbf{-0.1871 $\pm$ j0.7118}\\ -0.1978 $\pm$ j7.5172 \\ -0.2035 $\pm$ j6.0243 \\-0.2118 $\pm$ j5.4177}\\\hline 
10        &\tabincell{l}{\textbf{-0.1841 $\pm$ j6.4444}\\ -0.1955 $\pm$ j4.1341\\  -0.2075 $\pm$ j0.4542\\ -0.2080 $\pm$ j0.6288}
          &\tabincell{l}{\textbf{-0.1776 $\pm$ j6.1778}\\ -0.2059 $\pm$ j5.4403 \\ -0.2095 $\pm$ j7.5394 \\-0.2164 $\pm$ j3.8262}\\\hline           
\end{tabular}
}
\end{center}
\end{table}

\begin{table*}[!t]
\caption{Calculation Time of SQP-GS on 39-Bus and 118-Bus Systems}
\label{table_CPU}
\begin{center}
{
\begin{tabular}{c| c c c|c c c|c cc|c cc}\hline
\multirow{3}{*}{Step}  & \multicolumn{6}{c|}{New England 39-bus system}                        & \multicolumn{6}{c}{IEEE 118-bus system}         \\
                  & \multicolumn{3}{c|}{Sample size = 30} & \multicolumn{3}{c|}{Sample size = 60} & \multicolumn{3}{c|}{Sample size = 30} & \multicolumn{3}{c}{Sample size = 60} \\\cline{2-13}
                  &\tabincell{c}{CPU(s)\\/call} &\tabincell{c}{Calls\\/iter.} &\tabincell{c}{CPU(s)\\/iter.}   &\tabincell{c}{CPU(s)\\/call}&\tabincell{c}{Calls\\/iter.} &\tabincell{c}{CPU(s)\\/iter.}      &\tabincell{c}{CPU(s)\\/call} &\tabincell{c}{Calls\\/iter.} &\tabincell{c}{CPU(s)\\/iter.}    &\tabincell{c}{CPU(s)\\/call} &\tabincell{c}{Calls\\/iter.} &\tabincell{c}{CPU(s)\\/iter.}       \\\hline
 Eigen. Analysis   & 0.004   & 32  & 0.13              &0.004    & 62 &0.24                   &0.103    &32  &3.29              &0.108  &65  &7.05                   \\
 Sensitivities    & 0.005    & 31  & 0.14              &0.005     & 61 &0.27                  &0.028     &31 &0.86               &0.029 &61  &1.79                     \\
 QP               &0.12     &  1   &0.12                &0.14    & 1  &0.14                    &0.62      &1  &0.62              &0.72   &1 &0.72                    \\\hline
 Iter. Times     &\multicolumn{3}{c|}{68}       &\multicolumn{3}{c|}{62}       &\multicolumn{3}{c|}{13}    &\multicolumn{3}{c}{12}\\\hline
 Total CPU (s)        &\multicolumn{3}{c|}{26.72 (8.97\textsuperscript{*})}     &\multicolumn{3}{c|}{40.60 (9.54\textsuperscript{*})}     &\multicolumn{3}{c|}{60.33 (9.94\textsuperscript{*})}     &\multicolumn{3}{c}{114.93 (10.49\textsuperscript{*})}\\\hline    
 \multicolumn{13}{l}{\textsuperscript{*}\footnotesize{Estimated CPU time if ideal parallel techniques are used to sample the gradients}}\\
 \multicolumn{13}{l}{\footnotesize{\; Estimated CPU time =\big(Eigen. Analysis CPU(s)/call + Sensitivities CPU(s)/call + QP CPU(s)/call\big)$\times$Iter. Times + the other time.}}\\
  \multicolumn{13}{l}{\footnotesize{\; The other time includes L-BFGS time in Step 8 of Algorithm \ref{algo1}, linear search time in Step 20 of Algorithm \ref{algo1}, and data input and result output time. }}
\end{tabular}
}
\end{center}
\end{table*}

\begin{figure}[H]
    \centering
    \includegraphics[width = 0.48\textwidth]{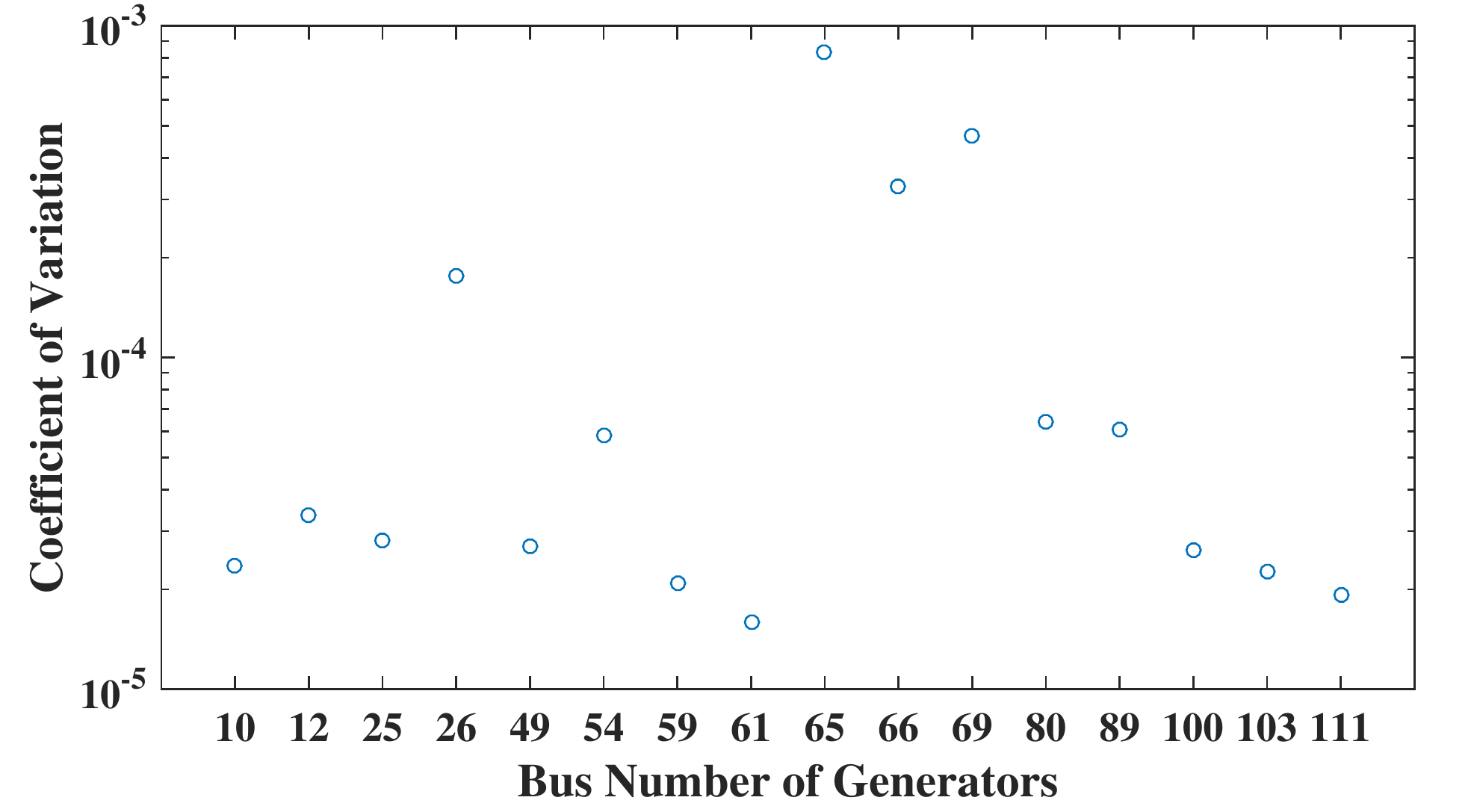}
    \caption{Coefficient of variation of generation outputs for 118-bus system over 20 runs of SSSC-OPF.}
    \label{fig:variance}
\end{figure}

Since the sampling in SQP-GS is random, the solutions may be different for different runs. 
We run SQP-GS for 20 times and the standard deviation of the generation cost is 0.06 $\$/$h. 
In Fig. \ref{fig:variance} we show the coefficient of variation, the standard deviation divided by the mean, of the generation output for each generator over 20 runs, which is  
very small and indicates that the difference between the solutions is small. 
\subsection{Efficiency}
We also test the efficiency of the proposed method on 39- and 118-bus systems. 
All simulations are carried out on a Dell Precision T5810  with a four-core 3.5 GHz processor and 64 GB RAM memory. 
Each iteration of the proposed method involves eigenvalue analysis, computing sensitivities, and solving QP subproblem.

Table \ref{table_CPU} lists the average step and the total CPU time with sample size $30$ and $60$. 
It is seen that the proposed SQP-GS method can solve the SSSC-OPF problem efficiently. 
Since in each iteration the gradient sampling performs eigenvalue analysis and sensitivity calculation for many times, these two steps are time consuming.
However, because the gradients can be sampled independently, the required CPU time can be greatly reduced by parallel computing. 
Furthermore, the computational burden of the eigenvalue analysis can be reduced by only calculating critical eigenvalues, such as by Jacobi-Davidson Method \cite{1583719}.

From Table \ref{table_CPU} it is also seen that the number of iterations for the IEEE 118-bus system is smaller than that for the New England 39-bus system.
The iteration times of SQP-GS can depend on the system size, the constraints, and the objective function.
In our test cases it seems that the small-signal stability constraint is the dominant factor, and 
the IEEE 118-bus case has a smaller number of iterations mainly because 
its $\overline{\eta}$ in the small-signal stability constraint is greater than that in the New England 39-bus case 
($\overline{\eta}=-0.1$ for IEEE 118-bus case and $\overline{\eta}=-0.2$ for New England 39-bus case).

\section{Conclusion} \label{conclusion}
In this paper, we propose an SQP-GS method to solve the SSSC-OPF problem, which is a nonsmooth optimization problem due to 
the property of the spectral abscissa function. 
In SQP-GS, a GS technique is used to evaluate the gradients around the current iterate to make the search
direction computation effective in the nonsmooth regions. 
The closed-form eigenvalue sensitivity is used to calculate the gradient of the spectral abscissa. 
Simulation results on three test systems show that the proposed method can solve the SSSC-OPF problem effectively 
without any convergence problem. 
By contrast, the existing Interior Point Method either cannot get as good solution as SQP-GS or even fail for relatively large systems.




%



\ifCLASSOPTIONcaptionsoff
  \newpage
\fi



\bibliographystyle{IEEEtran}
%
\bibliography{paperref}

%

\begin{IEEEbiography}[{\includegraphics[width=1in,height=1.25in,clip,keepaspectratio]{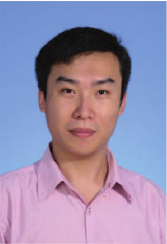}}]{Peijie Li}
(M'14) received the B.E. degree and Ph.D. degree in  electrical engineering from Guangxi University, Nanning, China, in 2006 and 2012, respectively.
   
From  2015, he is working at Argonne National Laboratory, Lemont, IL, USA, as a visiting scholar. He is currently also  an associate professor at the Guangxi University. His research interests include optimal power flow, small-signal stability, security constrained economic dispatch and restoration. 
\end{IEEEbiography}

\begin{IEEEbiography}[{\includegraphics[width=1in,height=1.25in,clip,keepaspectratio]{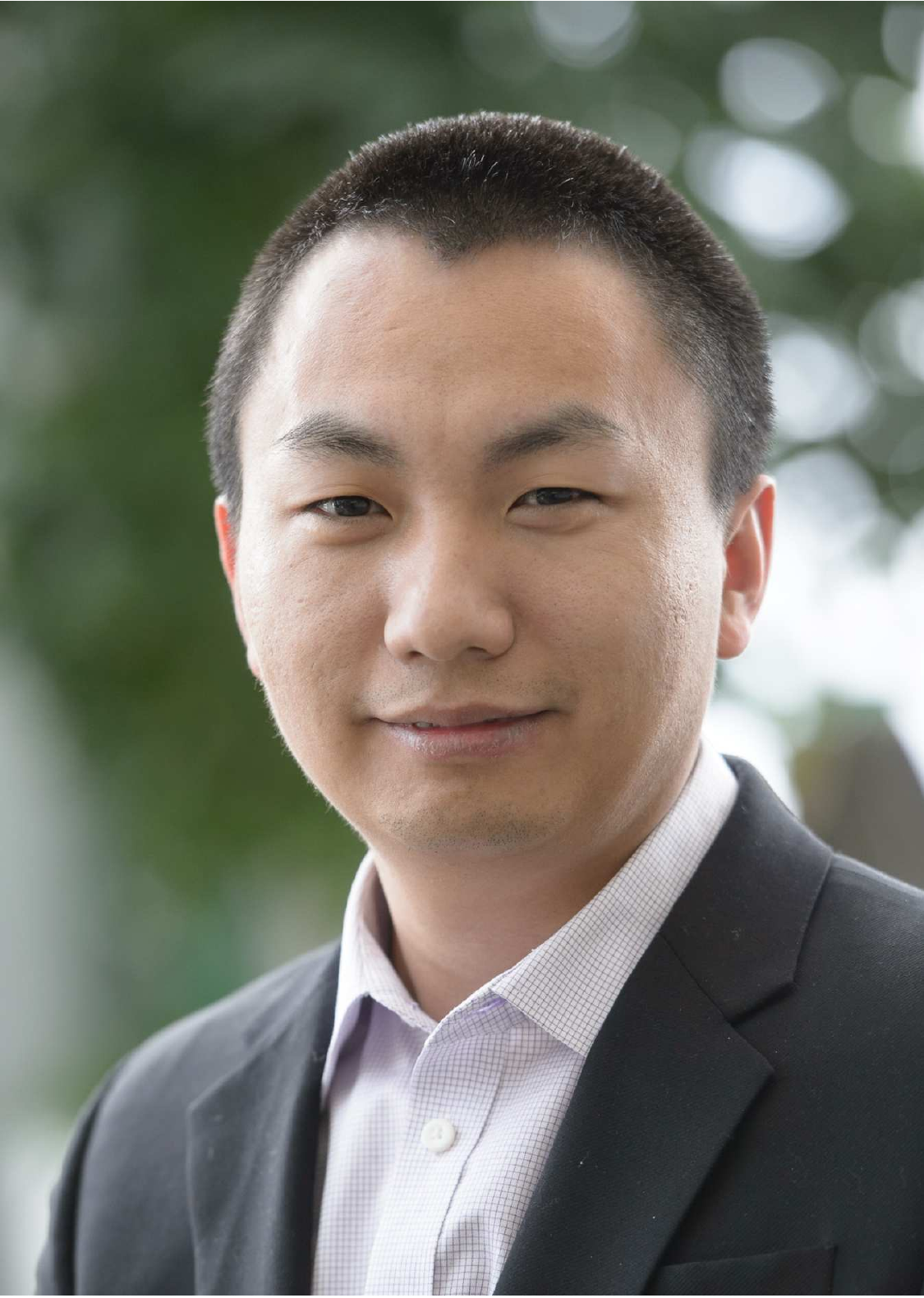}}]{Junjian Qi} (S'12--M'13)
received the B.E. degree from Shandong University, Jinan, China, in 2008 and the Ph.D. degree Tsinghua University, Beijing, China, in 2013, both in electrical engineering.
	
	In February--August 2012 he was a Visiting Scholar at Iowa State University, Ames, IA, USA. During September 2013--January 2015 he was 
	a Research Associate at Department of Electrical Engineering and Computer Science, University of Tennessee, Knoxville, TN, USA. 
	Currently he is a Postdoctoral Appointee at the Energy Systems Division, Argonne National Laboratory, Argonne, IL, USA. 
	His research interests include cascading blackouts, power system dynamics, state estimation, synchrophasors, and cybersecurity.
\end{IEEEbiography}

\begin{IEEEbiography}[{\includegraphics[width=1in,height=1.25in,clip,keepaspectratio]{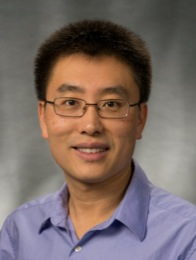}}]{Jianhui Wang}
(M'07-SM'12) received the Ph.D. degree in electrical engineering from Illinois Institute of Technology, Chicago, IL, USA, in 2007.

Presently, he is the Section Lead for Advanced Power Grid Modeling at the Energy Systems Division at Argonne National Laboratory, Argonne, IL, USA.

Dr. Wang is the secretary of the IEEE Power \& Energy Society (PES) Power System Operations Committee. He is an associate editor of Journal of Energy Engineering and an editorial board member of Applied Energy. He is also an affiliate professor at Auburn University and an adjunct professor at University of Notre Dame. He has held visiting positions in Europe, Australia and Hong Kong including a VELUX Visiting Professorship at the Technical University of Denmark (DTU). Dr. Wang is the Editor-in-Chief of the IEEE Transactions on Smart Grid and an IEEE PES Distinguished Lecturer. He is also the recipient of the IEEE PES Power System Operation Committee Prize Paper Award in 2015.
\end{IEEEbiography}

\begin{IEEEbiography}[{\includegraphics[width=1in,height=1.25in,clip,keepaspectratio]{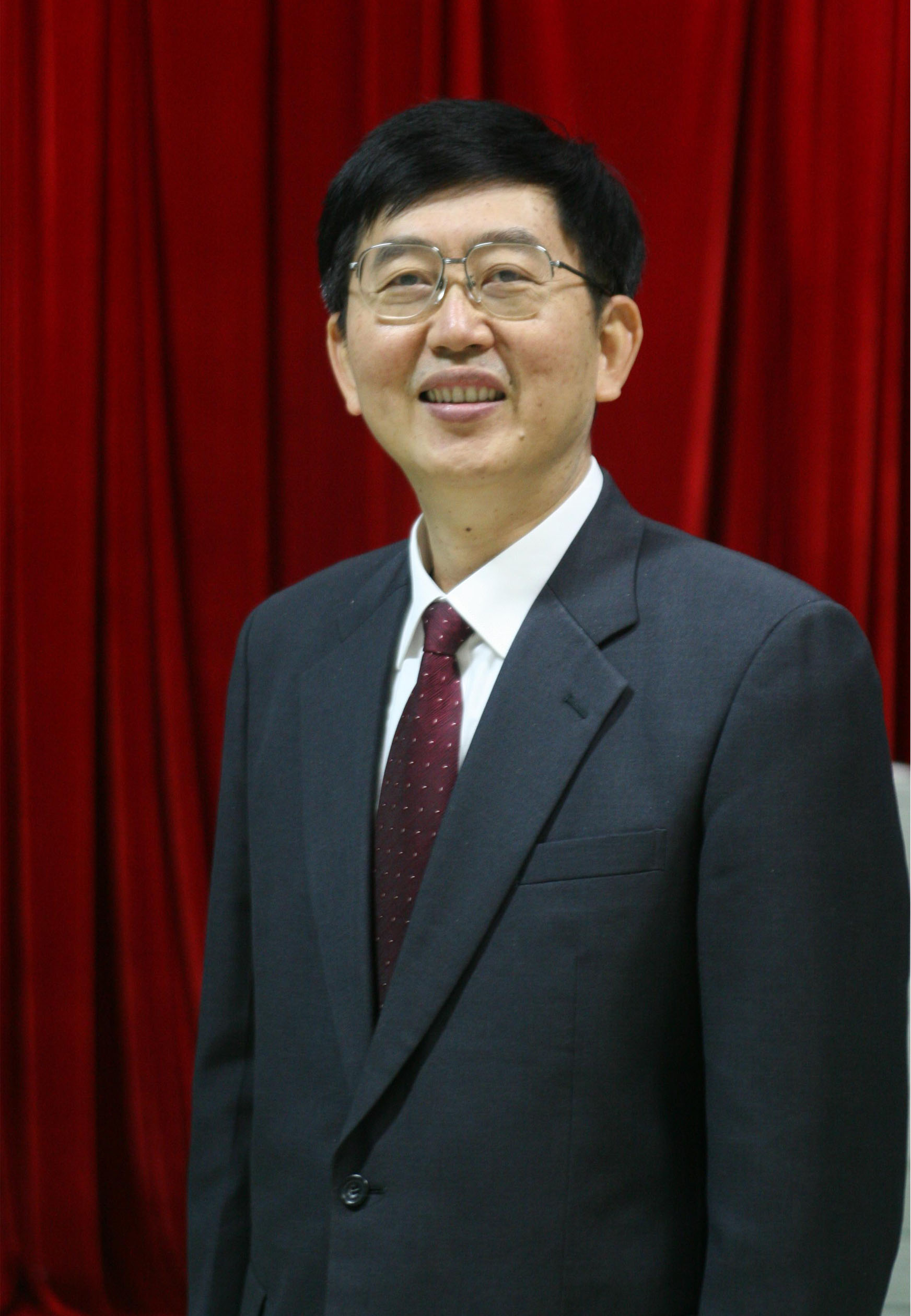}}]{Hua Wei}
  received BS and MS degrees in power
engineering from Guangxi university, Guangxi, P.R. China in 1981 and 1987, respectively, and Ph.D. degree in power
engineering from Hiroshima University, Japan in 2002. 

He was a visiting Professor at Hiroshima University, Japan from 1994 to 1997. From 2004 to 2014, he was the vice-president of the Guangxi University. Now, he is a  professor of Guangxi University. He is also the director  of the Institute of Power System Optimization, Guangxi University. His research interests include power system operation and planning, particularly in the application of optimization theory and methods to power systems.
\end{IEEEbiography}

\begin{IEEEbiography}[{\includegraphics[width=1in,height=1.25in,clip,keepaspectratio]{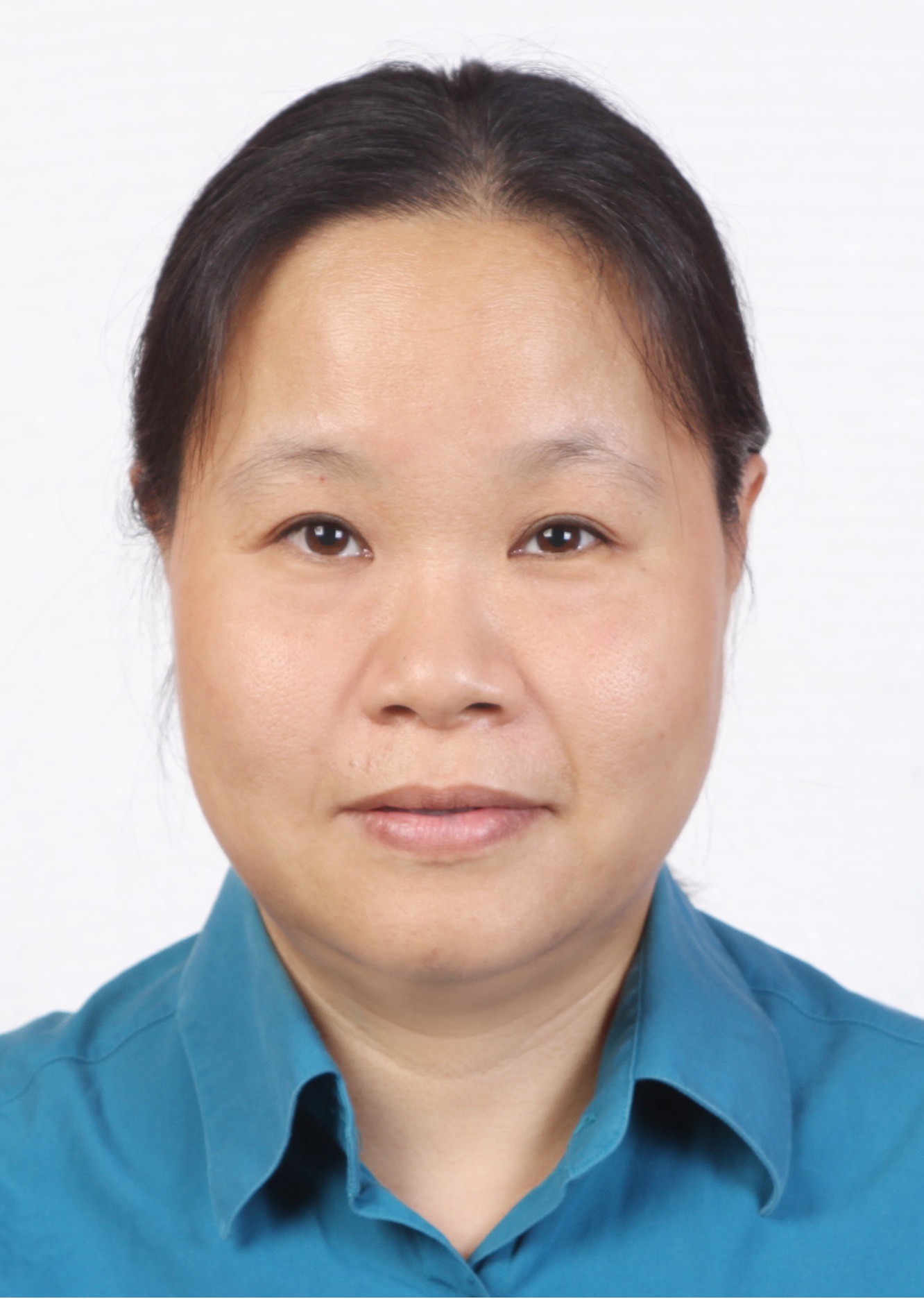}}]{Xiaoqing Bai}
 received B.S. and Ph.D. degrees in power engineering from Guangxi University,
Nanning, China, in 1991, 2010 respectively. 

She was a Post Doctoral Research Assistant of University of Nebraska-Lincoln from 2012 to 2015. Now she is a Professor of Institute  of Power System Optimization, Guangxi University. Her research interest is in power system optimization based on SDP and robust optimization.
\end{IEEEbiography}

\begin{IEEEbiography}[{\includegraphics[width=1in,height=1.25in,clip,keepaspectratio]{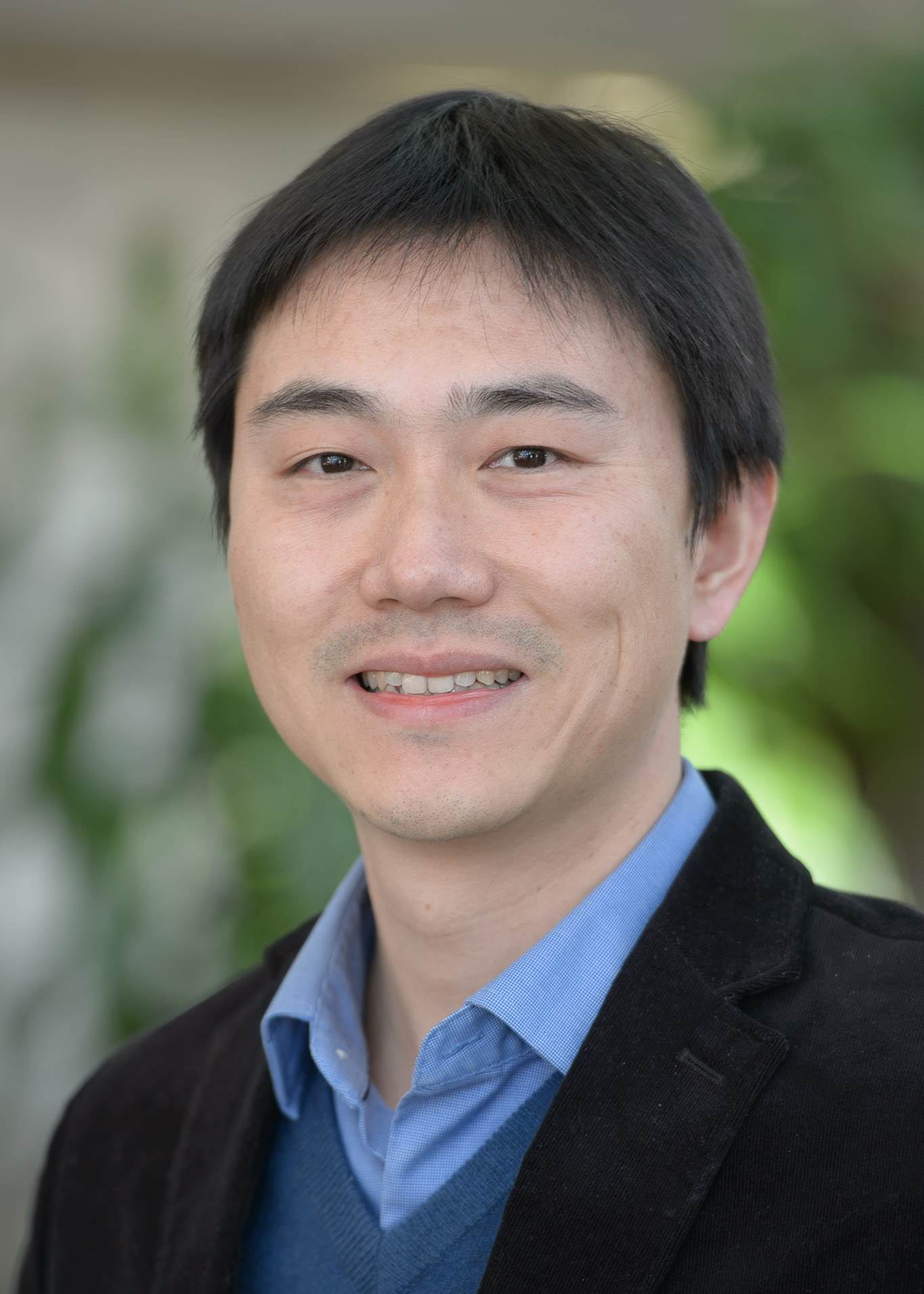}}]{Feng Qiu}
 (M'14) received his Ph.D. from the School of Industrial and Systems Engineering at the Georgia Institute of Technology in 2013. He is a computational scientist with the Energy Systems Division at Argonne National Laboratory, Argonne, IL, USA. His current research interests include optimization in power system operations, electricity markets, and power grid resilience. 
\end{IEEEbiography}





\end{document}